\documentclass[twoside,12pt]{article}
\usepackage{a4wide}
\usepackage{amsmath}
\usepackage{amsfonts}
\usepackage{amssymb}
\usepackage{latexsym}
\usepackage{times} 
\usepackage{pstricks}     

\def\date{25 July 2012}

\newcommand{\ed}{\end{document}}
\newcounter{mycnt}
\def\themycnt{\thesection.\arabic{mycnt}}
\def\mybenv#1{\refstepcounter{mycnt}%
       \vskip 3pt\noindent{\bf #1~~\themycnt}:~}
\def\myeenv{\hfill\rule{1ex}{1ex}\vskip 3pt}
\def\qed{\hfill$\Box$}
\def\ds{\displaystyle}
\def\ov{\overline}

 \def\nn{\nonumber \\}

\def\openR{\mathbb{R}}

\def\openN{\mathbb{N}}
\def\openZ{\mathbb{Z}}

\def\openQ{\mathbb{Q}}

\def\xb{\overline{x}}

\def\la{\langle\,}
\def\ra{\,\rangle}
\def\ot{\otimes}
\def\!{\kern -0.15ex}

\def\CR{\textsf{Char-}\!}
\def\SY{\textsf{Symm-}\!}

\def\grpGL{\mathsf{GL}}

\def\grpO{\mathsf{O}}
\def\grpSp{\mathsf{Sp}}
\def\grpSO{\mathsf{SO}}
\def\grpS{\mathcal{S}}
\def\grpG{\mathsf{G}}
\def\grpH{\mathsf{H}}
\def\grpHpi{\mathsf{H_{\pi}}}
\def\End{\text{\sf End}}

\def\conv{\star}

\def\ch{\text{ch}\,}

\def\[{[\,}
\def\]{\,]}

\def\A{\mathcal{A}}   
\def\B{\mathcal{B}}   %
\def\C{\mathcal{C}}   
\def\D{\mathcal{D}}   
\def\E{\mathcal{E}}   
\def\F{\mathcal{P}}   
\def\P{\mathcal{P}}

\def\mystrut{\hbox{\vrule height10.0pt depth2pt width0pt}}
\def\mybox{\hbox to 10.0pt}
\def\norulefill{\leaders\hrule height0pt\hfill}
\def\nr#1{\multispan{#1}\norulefill}
\def\hr#1{\multispan{#1}\hrulefill}

\def\mybigbox{\hbox to 15pt}

\def\mylargebox{\hbox to 20pt}

\def\K{K}  
\def\J{J}  
\def\grpGLrat{\mathsf{GLrat}}

\begin{document}
\title{\Large\sf The Hopf Algebra Structure of the Character Rings\\
of Classical Groups}
\author{Bertfried Fauser,~~Peter D. Jarvis~~and~~Ronald C. King}
{{\renewcommand{\thefootnote}{\fnsymbol{footnote}}
\footnotetext{\kern-15.3pt PACS numbers: 02.10.-v, 02.10.De, 02.20.-a, 02.20.Hj.
}
}}
\maketitle
\begin{abstract}
The character ring \CR$\grpGL$ of covariant irreducible tensor
representations of the general linear group admits a Hopf algebra structure isomorphic
to the Hopf algebra \SY$\Lambda$ of symmetric functions. Here we study the character 
rings \CR$\grpO$ and \CR$\grpSp$ of the orthogonal and symplectic subgroups of the 
general linear group within the same framework of symmetric functions. 
We show that \CR$\grpO$ and \CR$\grpSp$ 
also admit natural Hopf algebra structures that are isomorphic to that of \CR$\grpGL$,
and hence to \SY$\Lambda$. The isomorphisms are determined explicitly, along
with the specification of standard bases for \CR$\grpO$ and \CR$\grpSp$ analogous
to those used for \SY$\Lambda$. A major structural change arising from the adoption 
of these bases is the introduction of new orthogonal and symplectic Schur-Hall 
scalar products. Significantly, the adjoint with respect to multiplication 
no longer coincides, as it does in the \CR$\grpGL$ case, with a Foulkes derivative
or skew operation. The adjoint and Foulkes derivative now require separate definitions,
and their properties are explored here in the orthogonal and symplectic cases. 
Moreover, the Hopf algebras \CR$\grpO$ and \CR$\grpSp$ are not self-dual. The dual 
Hopf algebras \CR$\grpO^*$ and \CR$\grpSp^*$ are identified. Finally, the Hopf algebra 
of the universal rational character ring \CR$\grpGLrat$ of mixed irreducible tensor 
representations of the general linear group is introduced and its structure maps 
identified. 
\end{abstract}

\section{Introduction}
\label{Sec-intro}

\subsection{Motivation}
\label{Subsec-mot}

It is hardly possible to overestimate the importance of group representation
theory and the associated calculus of group characters. It plays a role in many
areas of physics, chemistry, biology and not least in pure mathematics. For that
reason, new techniques which deal with group characters in a unified and 
structural way are not only of interest in their own right, but also may be of
great help in more applied work.

A common problem involving the application of group representation theory is
that of determining an underlying \emph{symmetry group} whose irreducible
representations accommodate elementary particle or nucleon states. In this context 
\emph{tensor products} govern such things as interactions, scattering and decay 
processes. In addition \emph{symmetry breaking}, whereby the potential symmetries 
of some idealised system are more realistically limited to some subset of the original
symmetries, manifests itself group theoretically by way of \emph{restriction from 
group to subgroup} and the corresponding \emph{branching}, or \emph{reduction} of 
representations. Traditional homelands in physics for group theory-powered insights 
have included the multiplet organization of elementary particles, Wigner's 
nuclear multiplet theory, the nuclear interacting boson model, atomic and nuclear 
shell theory, as well as the building of grand unified theories. 

In all these cases, while the detailed construction of explicit group representations 
might be helpful, it is their characters that play the key role.
In the present paper we study group representations via the Hopf algebraic
structure of these characters. It is the \emph{products} and \emph{coproducts}
of these Hopf algebras that determine the decompositions of \emph{tensor products} 
and \emph{group-subgroup branching rules} that are required in physical 
applications, while the characters themselves, which as we shall see are all expressible
within the framework of the ring $\Lambda$ of symmetric functions, that specify the
physical states themselves. 

It was already observed in earlier work~\cite{fauser:jarvis:2003a,fauser:jarvis:king:wybourne:2005a} 
that Hopf algebra techniques allowed symmetric function methods to be organized and generalized in
an elegant way. The approach was developed in part by applying methods borrowed
from quantum field theory~\cite{fauser:2001b,brouder:fauser:frabetti:oeckl:2002a}, 
in a simplified group theoretical setting. In group theory terms, this earlier symmetric function
work concerns the characters of the general linear group.  In the present paper,
we pursue these investigations by turning to the \emph{classical subgroups} of
the general linear group. We show how the character rings of the orthogonal and
symplectic groups admit natural Hopf algebraic structures. We obtain these Hopf
algebras as isomorphic images of the Hopf algebra of the character ring of the
general linear group, which is in turn isomorphic to the Hopf algebra of
symmetric functions. The isomorphy is defined by the underlying branching, which
establishes an isomorphism between the module of characters of the general
linear group, and those of its classical subgroups. This module map induces a
map of Hopf algebras, as we are going to show.

Despite their isomorphism as Hopf algebras, the different character Hopf
algebras encode different information. This stems partly from the fact that
we are interested in \emph{canonical bases}, which \emph{differ} for 
different character modules. The prime example concerns the Schur functions
which furnish irreducible characters of the general linear group $\grpGL$. If
we branch from $\grpGL$ say to the orthogonal group $\grpO$, or the symplectic
group $\grpSp$, the Schur functions are no longer the irreducible characters,
and they lose, in part, their important and singular meaning. The orthogonality
of irreducible $\grpGL$ characters, corresponding to irreducible
representations with highest weight specified by integer partitions, $\lambda$, 
is expressed formally by means of the \emph{Schur-Hall scalar product} with respect to 
which the Schur functions $s_\lambda$ are orthonormal,
\begin{align}
\la s_\lambda \mid s_\mu\ra &= \delta_{\lambda,\mu}.
\end{align}
The decomposition of irreducible representations of $\grpGL$ on restriction
to the orthogonal or symplectic subgroups, $\grpO$ or $\grpSp$, involves a branching 
rule that is determined by expressing suitably restricted irreducible
characters of $\grpGL$ in terms of irreducible orthogonal or symplectic group
characters. These characters will be called Schur functions of
orthogonal or symplectic type. Since orthogonal and symplectic groups are
completely reducible, we can find a basis of such irreducible characters. It is hence
a group-theoretical necessity to introduce, on these character Hopf algebras,
\emph{new Schur-Hall scalar products} which express the fact that Schur
functions of orthogonal or symplectic type, $o_\lambda$ or $sp_\lambda$,
respectively, are mutually orthonormal (Schur's lemma):
\begin{align}
\la o_\lambda \mid o_\mu\ra_2 = \delta_{\lambda,\mu} &&\hbox{and}&&
\la sp_\lambda \mid sp_\mu\ra_{11} = \delta_{\lambda,\mu}\,.
\end{align}
The indexing stems from the plethystic origin of these particular branchings
(see below). These scalar products and the associated orthogonal bases
are the new structural elements which distinguish the otherwise isomorphic Hopf algebras.

The general case of symmetric function branchings was discussed 
in~\cite{fauser:jarvis:king:wybourne:2005a}. There we considered module
isomorphisms between the module of characters of a group $\grpG$ and the module
of characters of a subgroup $\grpH$. Specifically,  an algebraic subgroup
$\grpHpi$ of $\grpGL$ was taken, consisting of matrix  transformations fixing an
arbitrary tensor of symmetry type $\pi$ -- the orthogonal and  symplectic cases
correspond to the weight two symmetric and antisymmetric cases, $\pi=(2)$ and
$\pi=(1,1)$ respectively, of nonsingular bilinear forms.  However,
generically, the symmetric function bases obtained by branchings with respect
to higher rank invariants are no longer irreducible, but only indecomposable
at best. For this reason we study, in a first attempt, the orthogonal and
symplectic cases.

Even these classical cases reveal some novel features when treated in this
formal setting. We need to introduce new classes of Schur functions, as
described above, as is well-understood classically and was used at least
implicitly already by Weyl. Complete and elementary symmetric functions now have
different expansions in terms of irreducible orthogonal and symplectic characters; 
also it turns out that power sums pick up an extra additive term. More significantly, we
need to separate the notion of multiplicative adjoint (which we denote by
$s_\lambda^\dagger$), which leads to skew Schur functions in the $\grpGL$ case,
from that of the Foulkes derivative, which we denote by $s_\lambda^\perp$. This
stems from the fact that the adjoint of multiplication depends on the Schur-Hall
scalar product adopted, and that the branched Hopf algebras are no longer self
dual.

This work extends and elaborates on the material presented in the conference
paper~\cite{king:fauser:jarvis:2008a} with the addition of many proofs,
and thereby establishes the necessary tools
to deal, for example, with vertex operator algebras of orthogonal and symplectic type, as
described elswhere~\cite{fauser:jarvis:king:2010a}, see also for 
example~\cite{baker:1996a}. 
 
A further extension presented here covers rational characters of mixed tensor 
irreducible representations of the general linear group. In order to exploit Hopf algebra 
methods systematically in this context, and to make more rigorous previous discussions of 
products and branchings of mixed tensor representations~\cite{abramsky:king:1970a,king:1975a}, 
it is necessary to extend the underlying ring of symmetric functions from $\Lambda$ to 
$\Lambda\otimes\ov{\Lambda}$ and to define, following Koike~\cite{koike:1989a}, corresponding 
universal rational characters. 
 
\subsection{Organisation}
\label{Subsec-org}

The organisation of the paper is as follows. Some facts about the symmetric function 
Hopf algebra \SY$\Lambda$~\cite{fauser:jarvis:2003a} are provided in Section~\ref{Sec-symm-Hopf}, 
with an emphasis on its Schur function basis. This section also includes the definitions
of certain infinite series of Schur functions. These are used in Section~\ref{Sec-char-Hopf} 
to define the universal irreducible characters of the classical orthogonal and symplectic groups~\cite{littlewood:1940a,king:1975a,koike:terada:1987a} that are the main focus of our 
study. 

Thanks to their definition by way of certain branchings from the general linear 
group, whose universal irreducible characters are nothing other than Schur functions,
the corresponding orthogonal and symplectic character rings are actually Hopf algebras 
isomorphic to the universal Hopf algebra of symmetric functions 
\SY$\Lambda$~\cite{fauser:jarvis:king:wybourne:2005a}.
A complete directory is given in Table~\ref{tab1} of Section~\ref{Sec-char-Hopf} of the 
action of the structure maps
of these Hopf algebras as first described in abbreviated form in~\cite{king:fauser:jarvis:2008a}. 
Full proofs of all the structure map identities are provided here in Section~\ref{Sec-Hopf-proofs}. 
This is followed in Section~\ref{Sec-bases} by a discussion of the
\emph{different realizations} of the \emph{power sum}, \emph{complete} and 
\emph{elementary symmetric functions} in the orthogonal and symplectic cases,
noting where differences due to the underlying groups occur.  

The new orthogonal and symplectic scalar products are introduced in Section~\ref{Sec-scalar-dual},
in which we discuss the adjoint operation and the Foulkes derivative, and provide the correct
Hopf-theoretical definition of the latter, which allows applications to generic 
branchings. Finally, in this section, the dual Hopf algebras $\CR\grpO^*$ and $\CR\grpSp^*$ are
identified along with their structure maps. These dual Hopf algebras share the same structure
maps as the original Hopf algebras but are isomorphic to \SY$\Lambda$ extended so as to
include infinite series of Schur functions.

The universal rational characters associated with mixed tensor irreducible
representations of $\grpGL$ are introduced in Section~\ref{Sec-rational-universal}.
They form the basis of a ring $\Lambda\ot\ov{\Lambda}$ and an associated
new rational character Hopf algebra, $\CR\grpGLrat$. Explicit expressions
are derived for the corresponding products and coproducts, and the remaining Hopf
algebra structure maps are also identified explicitly.

\section{Symmetric functions and their Hopf algebra}
\label{Sec-symm-Hopf}

\subsection{The ring of symmetric functions $\Lambda$}
\label{Subsec-symm}

We use the standard notation of Macdonald's book~\cite{macdonald:1979a}. 
Symmetric functions are conveniently indexed by integer partitions 
$\lambda=(\lambda_1,\lambda_2,\ldots,\lambda_\ell,0,0,\ldots)$ with
$\lambda_1\geq\lambda_2\geq\cdots\geq\lambda_\ell>0$. The  $\lambda_i\in\openN$
are the parts of the partition, $\vert\lambda\vert = \sum_{i=1}^\ell
\lambda_i $ is the weight of the partition, while
$\ell(\lambda)=\ell$ is its length. Such a partition is often written 
without the trailing zeros. In exponent form 
$\lambda=(\ldots,k^{m_k},\ldots,2^{m_2},1^{m_1})$ where 
$m_k\geq0$ is the number of parts of $\lambda$ that are equal to $k$ for 
$k=1,2,\ldots$. With this notation it is convenient to
introduce $z_\lambda= \prod_{k\geq1}\, k^{m_k}\, m_k!$
  
Corresponding to each partition $\lambda$ there exists a Ferrers or Young diagram $F^\lambda$.
This consists of $|\lambda|$ boxes arranged in $\ell(\lambda)$ rows of 
lengths $\lambda_i$ for $i=1,2,\ldots,\ell(\lambda)$. 
The column lengths of $F^\lambda$ specify the
parts $\lambda'_j$ for $j=1,2,\ldots,\ell(\lambda')$
of the partition $\lambda'$ that is conjugate to $\lambda$.
If $F^\lambda$ has $r$ boxes on the
main diagonal, with arm lengths $\lambda_k-k=a_k$ and leg lengths $\lambda'_k-k=b_k$
for $k=1,2,\ldots,r$, then $\lambda$ is said to have rank $r(\lambda)=r$ and in Frobenius notation
$\lambda=\left(\begin{array}{llll}a_1&a_2&\cdots&a_r\\b_1&b_2&\cdots&b_r\end{array}\right)$
and 
$\lambda'=\left(\begin{array}{llll}b_1&b_2&\cdots&b_r\\a_1&a_2&\cdots&a_r\end{array}\right)$
with $a_1>a_2>\cdots>a_r\geq0$ and $b_1>b_2>\cdots>b_r\geq0$. 
Schematically,
we have
$$
{\vcenter
 {\offinterlineskip
 \halign{&\mystrut\vrule#&\mybox{\hss$\scriptstyle#$\hss}\cr
  \hr{15}\cr
                      & && && && && && && &\cr    
  \hr{15}\cr
                      & && && && &\cr       
  \hr{9}\cr
                      & && && && &\cr       
  \hr{9}\cr
                      & && &\cr          
  \hr{5}\cr
 }}}
=
{\vcenter
 {\offinterlineskip
 \halign{&\mystrut\vrule#&\mybox{\hss$\scriptstyle#$\hss}\cr
  \hr{15}\cr
                      &\lambda_1&\omit& &\omit& &\omit& &\omit& &\omit& &\omit& &\cr    
  \hr{15}\cr
                      &\lambda_2&\omit& &\omit& &\omit& &\cr       
  \hr{9}\cr
                      &\lambda_3&\omit& &\omit& &\omit& &\cr       
  \hr{9}\cr
                      &\lambda_4&\omit& &\cr          
  \hr{5}\cr
 }}}
=
{\vcenter
 {\offinterlineskip
 \halign{&\mystrut\vrule#&\mybox{\hss$\scriptstyle#$\hss}\cr
  \hr{15}\cr
   &\lambda_1^\prime&&\lambda_2^\prime&&\lambda_3^\prime&
   &\lambda_4^\prime&&\lambda_5^\prime& &\lambda_6^\prime& &\lambda_7^\prime&\cr    
  \hr{1}&\nr{7}&\hr{7}\cr
                                         & && && && &\cr    
  \hr{1}&\nr{7}&\hr{1}\cr
                                         & && && && &\cr    
  \hr{1}&\nr{3}&\hr{5}\cr
                                            & && &\cr    
  \hr{5}\cr
 }}}
=
{\vcenter
 {\offinterlineskip
 \halign{&\mystrut\vrule#&\mybox{\hss$\scriptstyle#$\hss}\cr
  \hr{15}\cr
                      & &&a_1&\omit& &\omit& &\omit& &\omit& &\omit& &\cr    
  \hr{15}\cr
                                 &b_1&& &&a_2&\omit& &\cr    
  \hr{1}&\nr{1}&\hr{7}\cr
                                      & &&b_2&& &&a_3&\cr    
  \hr{1}&\nr{3}&\hr{5}\cr
                                           & && &\cr    
  \hr{5}\cr
  }}}
$$
By way of an example we have
\begin{align}
 (7,4,4,2,0,\ldots) = (7,4^2,2)
 = \left(\begin{array}{rrr}
    6 & 2 & 1 \\
    3 & 2 & 0 
    \end{array}\right) \nonumber
\end{align}
and its conjugate
\begin{align}
(4,4,3,3,1,1,1,0\ldots)=(4^2,3^2,1^3)
 = \left(\begin{array}{rrr}
    3 & 2 & 0 \\
    6 & 2 & 1 
    \end{array}\right)\,.\nonumber    
\end{align}

Partitions are used to specify a number of objects of interest in the  present
work. Amongst these are the Schur functions $s_\lambda$. These form an
orthonormal $\openZ$-basis for the ring $\Lambda$ of symmetric functions.  To
be more precise, let $\openZ[x_1,\ldots,x_N]$ be the polynomial ring, or the
ring of formal power series, in $N$ commuting variables $x_1,\ldots,x_N$. The symmetric group
$\grpS_N$ acting on  $N$ letters acts on this ring by permuting the variables. 
For $\pi\in \grpS_N$ and $f \in \openZ[x_1,\ldots,x_N]$ we have
\begin{align}
\pi f(x_1,\ldots,x_N) &= f(x_{\pi(1)},\ldots,x_{\pi(N)}) \,.
\end{align}
We are interested in the \emph{subring of functions} invariant under this
action, $\pi f = f$, that is to say the ring of symmetric polynomials in
$N$ variables: 
\begin{align}
\Lambda_N[x_1,\ldots,x_N]&= \openZ[x_1,\ldots,x_N]^{S_N}.
\end{align}
This ring may be graded by the degree of the polynomials, so that
\begin{align}
\Lambda_N[x_1,\ldots,x_N]=\oplus_n\ \Lambda_N^{(n)}[x_1,\ldots,x_N] \,,
\end{align} 
where $\Lambda_N^{(n)}[x_1,\ldots,x_N]$ consists of homogenous symmetric
polynomials  in $x_1,\ldots,x_N$ of total degree $n$. 

In order to work with an arbitrary number of variables, following
Macdonald~\cite{macdonald:1979a}, we define the ring of symmetric functions
$\Lambda=\lim_{N\rightarrow\infty}\Lambda_N$ in its stable limit  ($N\rightarrow
\infty$) where  $\Lambda_N=\Lambda_M[x_1,\ldots,x_N,0,\ldots,0]$ for all $M\geq
N$. This ring of symmetric functions inherits the grading $\Lambda=\oplus_n\,
\Lambda^{(n)}$, with $\Lambda^{(n)}$ consisting of homogeneous symmetric
polynomials  of degree $n$.

A $\openZ$ basis of $\Lambda^{(n)}$ is provided by the monomial symmetric
functions $m_\lambda$ where $\lambda$ is any partition of $n$. There exist
further (integral and rational) bases for $\Lambda^{(n)}$ that are indexed by
the partitions $\lambda$ of $n$.  These are the complete, elementary and power
sum symmetric function bases defined \emph{multiplicatively} in terms of
corresponding one part functions by
\begin{align}
h_\lambda&=h_{\lambda_1}h_{\lambda_2}\ldots h_{\lambda_l}\,,
&&&
e_\lambda&=e_{\lambda_1}e_{\lambda_2}\ldots e_{\lambda_l}\,,
&&&
p_\lambda&=p_{\lambda_1}p_{\lambda_2}\ldots p_{\lambda_l},
\label{Eq-hep}
\end{align}
where the one part functions are defined for all $n\in\openN$ by 
\begin{align}
h_n&=\sum_{i_1\leq i_2\leq \cdots \leq i_n}\ x_{i_1} x_{i_2} \cdots x_{i_n}\,,
&&&
e_n&=\sum_{i_1 < i_2 < \cdots < i_n}\ x_{i_1} x_{i_2} \cdots x_{i_n}\,,
&&&
p_n&=\sum_i x_i^n.
\label{Eq-hnenpn}
\end{align}
With the convention $h_0=e_0=p_0=1$ and $h_{-n}=e_{-n}=p_{-n}=0$ for
positive $n$, their generating functions take the form
\begin{align}
H_t \ =\ \prod_{i\ge 1} \frac{1}{(1-x_it)}  &= \sum_{n\ge0}h_n t^n, &&& 
E_t \ =\ \prod_{i\ge 1} (1+x_it) &= \sum_{n\ge 0} e_n t^n, &&& 
t \frac{\mathrm{d}}{\mathrm{d}t} \log H_t & = \sum_{n\ge 1} p_n t^{n}. 
\label{Eq-hepn}
\end{align}

The most important \emph{non-multiplicative} basis of $\Lambda^{(n)}$ is provided
by the Schur functions $s_\lambda$ with $\lambda$ running over all partitions
of $n$. For a finite number of variables the Schur function
$s_\lambda(x_1,\ldots,x_N)$ may be defined as a ratio of alternants. It is a
homogeneous symmetric polynomial of total degree $n$, and is stable in the sense that 
$s_\lambda(x_1,\ldots,x_N,0,\ldots,0)=s_\lambda(x_1,\ldots,x_N)$ regardless of
how many $0$'s are appended to the list of variables. Taking the limit as
$N\rightarrow\infty$ of $s_\lambda(x_1,\ldots,x_N)$  serves to define the
required $s_\lambda\in\Lambda^{(n)}$~\cite{macdonald:1979a}. 

Varying $\lambda$ over all partitions, the Schur functions $s_\lambda$ provide a
$\openZ$-basis of $\Lambda$. We can go further. There exists a bilinear form on
$\Lambda$, the Schur-Hall scalar product $\la\cdot\mid\cdot\ra$. With respect to
this scalar product,  the Schur functions form an orthonormal basis of
$\Lambda$. In fact we have:
\begin{align}
\la s_\lambda\mid s_\mu \ra 
=\delta_{\lambda,\mu}\,, \qquad 
 \la p_\lambda\mid p_\mu\ra 
 & = z_\lambda \delta_{\lambda,\mu}\,, \qquad
\la m_\lambda\mid h_\mu \ra
&=\delta_{\lambda,\mu}\,, \qquad
\la f_\lambda\mid e_\mu \ra
&=\delta_{\lambda,\mu}. 
\end{align} 
These relations serve to define the `monomial' symmetric functions $m_\lambda$,
and the so-called `forgotten' symmetric functions $f_\lambda$ (for details 
see~\cite{macdonald:1979a}). 

In what follows we make use of various notation for Schur functions, including for example 
$s_\lambda(x_1,\ldots,x_N)$, $s_\lambda(x)$, $s_\lambda(x,y)$ or $s_\lambda$, 
depending on whether or not it is necessary to be explicit about the number
of  variables or the sets of variables under consideration. Here, a single
symbol $x$  may often stand for an alphabet, $x_1,x_2,\ldots$, finite or
otherwise,  while a pair $x,y$ signifies a pair of such alphabets
$x_1,x_2,\ldots,y_1,y_2,\ldots\,$.

\subsection{The Hopf algebra $\SY\Lambda$}
\label{Subsec-Hopf}

The graded ring of symmetric functions $\Lambda$ spanned by the Schur functions
$s_\lambda$ affords a graded self-dual, bicommutative Hopf algebra, which we
denote by $\SY\Lambda$, as can be seen once we have identified the appropriate product,
coproduct, unit, counit, antipode and self-duality condition. This can be done
as follows.

The {\bf outer product} of Schur functions is given in prefix form, 
infix dot product form, first without and then with variables, and 
finally in explicit form:
\begin{align}
m(s_\mu \ot s_\nu)&=s_\mu \cdot s_\nu \,, 
\nn
m(s_\mu(x)\ot s_\nu(y))
&=s_\mu(x)\cdot s_\nu(x) 
=\sum_{\lambda} C^\lambda_{\mu,\nu} s_\lambda(x)\,.
\label{Eq-dot-prod}
\end{align}
In some situations the $\cdot$ is omitted in favour of simple juxtaposition. 
Here, and elsewhere if not otherwise specified, tensor products are over $\openZ$ 
(or $\openQ$ if power sums are involved).

The {\bf outer coproduct} map is denoted by $\Delta$, and we use the variable
or the tensor product notation interchangeably, to give the coproduct in
prefix form, Sweedler index form~\cite{sweedler:1969a} and skew product forms, 
first without and then with variables, and finally in explicit form:
\begin{align}
\Delta(s_\lambda) 
&= s_{\lambda_{(1)}}\ot s_{\lambda_{(2)}}
=\sum_\nu s_{\lambda/\nu}\ot s_\nu=\sum_\mu s_{\mu}\ot s_{\lambda/\mu} \,,\nn
\Delta(s_\lambda(x))
&= s_\lambda(x,y) = s_{\lambda_{(1)}}(x)s_{\lambda_{(2)}}(y) 
= \sum_{\mu,\nu} C^{\mu,\nu}_\lambda s_\mu(x) s_\nu(y)\,.
\label{Eq-skew_prod}
\end{align}

The coefficient $C^\lambda_{\mu,\nu}$ of the multiplication map $m$ in the 
Schur function basis, and the structure constant $C_\lambda^{\mu,\nu}$ of  
the coproduct in the same basis turn out to be identical. 
This equality of coefficients
is a consequence of the self-duality condition
\begin{align}
  \la s_\lambda\mid m(s_\mu\ot s_\nu)\ra &=
  \la\Delta(s_\lambda)\mid s_\mu\ot s_\nu\ra \,. 
\label{Eq-duality}
\end{align}
In fact, although logically distinct, they are both equal to the famous
Littlewood-Richardson coefficient
$c^\lambda_{\mu,\nu}$, which may be evaluated combinatorially using
the Littlewood-Richardson rule~\cite{littlewood:1940a,macdonald:1979a}.
Thus we have
\begin{align} 
C_{\mu,\nu}^\lambda=c_{\mu,\nu}^\lambda=C^{\mu,\nu}_\lambda\,.
\end{align}
This follows from the well know fact that the dot and skew products 
of Schur functions are adjoint with respect to the Schur-Hall scalar product, 
that is to say~\cite{macdonald:1979a} 
\begin{align}
   \la s_\lambda\mid s_\mu\cdot s_\nu\ra = c^\lambda_{\mu,\nu}=\la s_{\lambda/\nu}\mid s_\mu\ra\,.
\label{Eq-dot-skew}   
\end{align}
Here, in the Schur function basis the operations of outer multiplication and that of skewing
are both defined in terms of Littlewood-Richardson coefficients by
\begin{align}
s_{\mu}\cdot s_{\nu}  
=\sum_\lambda\ c^\lambda_{\mu,\nu}\ s_\lambda
\qquad \mbox{and} \qquad
s_\mu^\perp(s_\lambda)=s_{\lambda/\mu}= \sum_\nu\ c^\lambda_{\mu,\nu}\ s_\nu\,.
\end{align}
The notation $s_\mu^\perp(s_\lambda)$ has been introduced to emphasise 
the fact that the ring of symmetric functions has a module structure 
under ${}^\perp$:
\begin{align}
f^\perp(g^\perp (h)) = (gf)^\perp(h) \qquad \mbox{or equivalently} \qquad
(h/g)/f = h/(gf).
\end{align}

The {\bf unit map} $\eta$, {\bf counit map} $\epsilon$, and {\bf antipode} 
$S$, are defined by:
\begin{align}
\eta:1\rightarrow s_0\,,\qquad
\epsilon:s_\lambda\rightarrow\delta_{\lambda,0}\,,\qquad
S:s_\lambda\rightarrow (-1)^{|\lambda|}s_\lambda'\,.
\end{align}
It is important to note that the following antipode identity in the Hopf
algebra \SY$\Lambda$:
\begin{align}
m\,(I\ot S)\Delta\,(s_\lambda)&=\eta\,\epsilon(s_\lambda)
\label{Eq-mIS}
\end{align}
yields the result~\cite{fauser:jarvis:2003a}
\begin{align}
\sum_\nu\ (-1)^{|\nu|}\,s_{\lambda/\nu}\cdot s_{\nu'}
&=\delta_{\lambda,0}\, s_0\,,
\label{Eq-antipode-identity}
\end{align}
since
\begin{align}
m(I\ot S)\Delta(s_\lambda)&=m(I\ot S)(\sum_\nu\ s_{\lambda/\nu}\ot s_\nu)\nn
&=m(\sum_\nu\ s_{\lambda/\nu}\ot(-1)^{|\nu|}s_{\nu'})
=\sum_\nu\ (-1)^{|\nu|}\,s_{\lambda/\nu}\cdot s_{\nu'}
\nonumber
\end{align}
and
\begin{align}
\eta\,\epsilon(s_\lambda)&=\eta(\delta_{\lambda,0})=\delta_{\lambda,0}\, s_0\,.
\nonumber
\end{align}

Returning to the bases provided by $h_\lambda$, $e_\lambda$ and $p_\lambda$
in (\ref{Eq-hep}), these bases are so-called \emph{multiplicative}, because
the outer product is just the \emph{(unordered) concatenation product}. Using 
self-duality, this means that the coproduct is just deconcatenation of these
products, together with the action:
\begin{align}
\Delta(h_n)&=\sum_{r=0}^n h_{n-r}\ot h_r\,, &&&
\Delta(e_n)&=\sum_{r=0}^n e_{n-r}\ot e_r\,, &&&
\Delta(p_n)&=p_n\ot 1 + 1 \ot p_n\,.
\label{Eq-coprod-hep}
\end{align}
These results all follow immediately from the definitions of (\ref{Eq-hepn}).
The first two of these show that one part complete and elementary
symmetric functions are divided powers. The third shows that the one part
power sum symmetric functions are the primitive elements of the Hopf
algebra \SY$\Lambda$.

\subsection{Schur function series}
\label{Subsec-series}

To describe characters of the orthogonal and symplectic groups effectively, 
Littlewood~\cite{littlewood:1940a}  introduced a set of infinite series of
Schur functions which we are frequently going to use; 
consult also~\cite{king:1975a,black:king:wybourne:1983a}. 
It is convenient to extend our ring $\Lambda$
to $\Lambda[[t]]$, where $t$ is a formal parameter, and our Hopf
algebra to $\SY\Lambda[[t]]$ itself extended so as to encompass infinite
series of Schur functions. Some of these Schur function series read
\begin{align}
A_t &=\sum_{\alpha\in\A} (-1)^{\vert\alpha\vert/2} t^{|\alpha|}\, \{\alpha\}
  &&&
B_t &=\sum_{\beta\in\B}  t^{|\beta|}\, \{\beta\}
\nn
C_t &=\sum_{\gamma\in\C} (-1)^{\vert\gamma\vert/2} t^{|\gamma|}\, \{\gamma\}
  &&&
D_t &=\sum_{\delta\in\D} t^{|\delta|}\,  \{\delta\}
\nn
E_t &=\sum_{\epsilon\in\E} (-1)^{(\vert\epsilon\vert+r(\epsilon))/2} t^{|\epsilon|}\, \{\epsilon\}
  &&&  
F_t &= \sum_{\zeta\in\F} t^{|\zeta|}\,  \{\zeta\}  
  &&&
\nn
G_t &=\sum_{\epsilon\in\E} (-1)^{(\vert\epsilon\vert-r(\epsilon))/2} t^{|\epsilon|}\, \{\epsilon\}
  &&&
H_t &=\sum_{\zeta\in\F} (-1)^{\vert\zeta\vert} t^{|\zeta|}\, \{\zeta\}  
\nn
L_t &=\sum_{m\geq0} (-1)^m t^m\, \{1^m\}
  &&&
M_t &=\sum_{m\geq0} t^m\,\{m\} 
  \nn
P_t &=\sum_{m\geq0} (-1)^m t^m\,\{m\}
  &&&
Q_t &=\sum_{m\geq0} t^m\,\{1^m\} 
\label{Eq-series}
\end{align}
where $m$ is summed over all non-negative integers, 
while $\P$ is the set of all partitions, $\D=2\P$ is the 
set of partitions all of whose parts are even, and
$\B$ is the set of partitions that are conjugate to those of $\D$. To define
$\A$, $\C$ and $\E$ it is convenient   
in Frobenius notation to let 
\begin{align}
\mathcal{P}_n\ = \left\{
\left(\begin{array}{cccc}
    a_1 & a_2 & \ldots & a_r \\
    b_1 & b_2 & \ldots & b_r 
    \end{array}\right) \in \mathcal{P}\ \,\bigg|\,   a_k-b_k=n  
    \quad \hbox{for all}\quad 
    \begin{array}{l} r=0,1,2,\ldots \\ k=1,2,\ldots,r\\ \end{array} \right\}
\label{Eq-frobenius-sets}
\end{align}
for all integers $n$. With this notation $\A=\P_{-1}$, $\C=\P_1$ and $\E=\P_0$.
Thus the partitions in $\C$ are the conjugates of those in $\A$, while
$\E$ is the set of all self-conjugate partitions. 
It should be pointed out that each of the Schur function series of (\ref{Eq-series})
includes the term $\{0\}=s_0=1$ since $\{1^0\}=\{0\}$ and 
the sets $\P$, $2\P$ and $\P_n$, for all integers $n$, contain the zero partition $(0)$ 
for which $|(0)|=\ell(0)=r(0)=0$ and $(0)'=(0)$.

The generating functions which serve to define these series take the form:
\begin{align}
A_t &:=\prod_{i<j} (1-t^2x_ix_j)
  &&&
B_t &:=\prod_{i<j} (1-t^2x_ix_j)^{-1}
\nn
C_t &:=\prod_{i\leq j} (1-t^2x_ix_j)
  &&&
D_t &:=\prod_{i\leq j} (1-t^2x_ix_j)^{-1}
\nn
E_t &:=\prod_k (1-t\,x_k)\prod_{i<j} (1-t^2x_ix_j)
  &&&  
F_t &:=\prod_k (1-t\,x_k)^{-1}\prod_{i<j} (1-t^2x_ix_j)^{-1} 
  &&&
\nn
G_t &:=\prod_k (1+t\,x_k)\prod_{i<j}(1-t^2x_ix_j)
  &&&
H_t &:=\prod_k (1+t\,x_k)^{-1}\prod_{i<j}(1-t^2x_ix_j)^{-1}
\nn
L_t &:=\prod_k (1-t\,x_k) 
  &&&
M_t &:=\prod_k (1-t\,x_k)^{-1} 
\nn
P_t &:=\prod_k (1+t\,x_k)^{-1}  
  &&&
Q_t &:=\prod_k (1+t\,x_k) 
\label{Eq-series-genfun}
\end{align}
As can be seen there is some redundancy here because $P_t=M_{-t}$, $Q_t=L_{-t}$,
$G_t=E_{-t}$ and $H_t=F_{-t}$, however we keep all 12 series $Z_t$ because in 
what follows we often denote $Z_1$ by $Z$. On the other hand we sometimes need
to display the arguments $x=(x_1,x_2,\ldots)$ of the above series by adopting
the more explicit notation $Z_t(x)$. 

A major feature of the
above list of Schur function series is that they come, as can be seen from their 
generating functions, in mutually inverse pairs:
\begin{align}
A_tB_t=C_tD_t=E_tF_t=G_tH_t=L_tM_t=P_tQ_t=1.
\label{Eq-series-pairs}
\end{align}
Moroever
\begin{align}
L_tA_t=P_tC_t=E_t,~~~M_tB_t=Q_tD_t=F_t,~~~M_tC_t=Q_tA_t=G_t,~~~L_tD_t=P_tB_t=H_t.
\label{Eq-series-identities}
\end{align}

The above generating functions make it particularly easy to establish the following:

\mybenv{Proposition}
Let $Z_t$ be any one of the Schur function series (\ref{Eq-series-genfun}), then
\begin{align}
  \Delta(Z_t) = (Z_t\otimes Z_t)\ \Delta'(Z_t) 
\end{align}
where the cut coproducts are given by
\begin{align}
    \Delta'(Z_t) &= \left\{  
    \begin{array}{ccl} 
   \ds\sum_{\sigma\in\P} (-t^{2})^{|\sigma|}\, \{\sigma\} \otimes \{\sigma'\} 
   &\hbox{for}& Z_t=A_t,C_t,E_t,G_t;\cr\cr
   \ds \sum_{\sigma\in\P} \, t^{2\,|\sigma|}\, \{\sigma\} \otimes \{\sigma\} 
   &\hbox{for}& Z_t=B_t,D_t,F_t,H_t;\cr\cr
    1&\hbox{for}& Z_t=L_t,M_t,P_t,Q_t.
    \end{array} \right.
\end{align}
\label{Prop-series-coproducts}
\myeenv

\paragraph{Proof:}
Let $z=(z_1,z_2,\ldots)=(x,y)=(x_1,x_2,\ldots,y_1,y_2,\ldots)$
and note that $\Delta(Z_t(z))=Z_t(x,y)=Z_t(x)Z_t(y)Z_t'(x,y)$ where the generating functions
immediately imply that for $Z_t=A_t,C_t,E_t,G_t$ we have $Z_t'(x,y)=\prod_{i,j} (1-t^2x_iy_j)$,
and for $Z_t=B_t,D_t,F_t,H_t$ we have $Z_t'(x,y)=\prod_{i,j} (1-t^2x_iy_j)^{-1}$, whilst
for $Z_t=L_t,M_t,P_t,Q_t$ we have $Z_t'(x,y)=1$. It only remains to use the well known 
Cauchy identity~\cite{macdonald:1979a}
\begin{align}
\label{Eq-Cauchy}
 \K_t(x,y):=\prod_{i,j} (1-t^2x_iy_j)^{-1} = \sum_{\sigma\in\P} \,t^{2\,|\sigma|} \{\sigma\} \otimes \{\sigma\}
\end{align}
and its dual
\begin{align}
\label{Eq-dual-Cauchy}
 \J_t(x,y):=\prod_{i,j} (1-t^2x_iy_j) = \sum_{\sigma\in\P} (-t^{2})^{|\sigma|}\, \{\sigma\} \otimes \{\sigma'\}\,.
\end{align}
\qed

It might be added here that
\begin{align}
    \prod_{i,j} (1-t^2 x_ix_j)^{-1} = B_t(x)\ D_t(x) 
    = \sum_{\sigma\in\P} \,t^{2\,|\sigma|} s_\sigma(x) s_{\sigma}(x)
    = \sum_{\sigma\in\P} \,t^{2\,|\sigma|} \{\sigma\cdot\sigma\}(x)
\end{align}
and
\begin{align}
    \prod_{i,j} (1-t^2 x_ix_j) = A_t(x)\ C_t(x) 
    = \sum_{\sigma\in\P} (-t^{2})^{|\sigma|} s_\sigma(x) s_{\sigma'}(x)
    = \sum_{\sigma\in\P} (-t^{2})^{|\sigma|} \{\sigma\cdot\sigma'\}(x)\,.
\end{align}

We are now in a position to exploit these series and their associated identities
in the specification of characters of the classical groups, in particular
what are known as their universal characters, and to study the detailed 
properties of the Hopf algebras of their character rings, initiated 
in~\cite{fauser:jarvis:king:wybourne:2005a} and~\cite{king:fauser:jarvis:2008a}.

\section{The Hopf algebras of the character rings of classical groups}
\label{Sec-char-Hopf}

\subsection{Universal characters of covariant tensor representations}
\label{Subsec-univ}

The groups, $\grpG$, under consideration here are the general linear group
$\grpGL$, the orthogonal group $\grpO$ and the symplectic group $\grpSp$.  If
the classical groups $\grpGL$, $\grpO$ and $\grpSp$ act by way of linear 
transformations in a space $V$ of dimension $N$, then they are denoted by
$\grpGL(N)$, $\grpO(N)$ and $\grpSp(N)$, respectively. Initially we confine attention to
their finite-dimensional irreducible covariant tensor  representations
$V_\grpG^\lambda$. Each of these is specified by their \emph{highest weight}
$\lambda$, which in each case is a partition. The corresponding character is 
denoted by $\ch V_\grpG^\lambda$. These characters may be expressed  by
means of Weyl's character formula~\cite{weyl:1939a} in terms of the eigenvalues
$(x_1,\ldots,x_N)$ of each group element $g\in\grpG$ realised as a matrix 
$M\in \End(V)$ of linear transformations in $V$.

It is well known that in the case of $\grpGL(N)$ we have 
\begin{align}
   \ch V_{\grpGL(N)}^\lambda &=\{\lambda\}(x_1,\ldots,x_N)=s_\lambda(x_1,\ldots,x_N)\,,
\label{Eq-charG}
\end{align}
where the central symbol accords with the notation of
Littlewood~\cite{littlewood:1940a}.  This character shares the same stable
$N\rightarrow\infty$ limit as Schur functions, and in this limit we define the
\emph{universal character}~\cite{koike:terada:1987a,king:1989a}
\begin{align}
   \ch V_{\grpGL}^\lambda &=\{\lambda\}=s_\lambda\,.
\label{Eq-charGL}
\end{align}

The orthogonal and symplectic groups leave invariant a symmetric second rank
tensor $g_{ij}=g_{ji}$ and an antisymmetric second rank tensor $f_{ij}=-f_{ji}$,
respectively. It is necessary to distinguish between the even and odd cases:
$N=2K$ and $N=2K+1$ with $K\in\openN$. The groups $\grpO(2K)$, $\grpO(2K+1)$ and
$\grpSp(2K)$ are all reductive Lie groups whose finite-dimensional
representations are fully reducible. On the other hand $\grpSp(2K+1)$, an
odd-dimensional symplectic group, is not reductive. This is a consequence of
the fact that its invariant bilinear form is singular. It can be realised as an
affine extension of $\grpSp(2K)\times \grpGL(1)$, that is the semi-direct
product of these groups with a set of translations as explained  by
Proctor~\cite{proctor:1988a}. As a result its finite-dimensional representations
are not necessarily fully reducible. Indeed its defining representation, $V$, of 
dimension $2K+1$ is indecomposable but contains two irreducible constituents of
dimensions $2K$ and $1$. More generally, Proctor has established that the
representations $V_{\grpSp(2K+1)}^\lambda$ are reducible but indecomposable  for
$\lambda\neq0$.

Despite these issues associated with the evenness and oddness of $N$, there still
exists a stable $N\rightarrow\infty$ limit and associated universal 
characters~\cite{koike:terada:1987a,king:1989a} denoted here by
$\ch V_{\grpO}^\lambda$ and $\ch V_{\grpSp}^\lambda$.
The Schur function series that we have introduced enable us to write
down Schur function expressions for these universal characters of $\grpO$
and $\grpSp$ in the form~\cite{littlewood:1940a,king:1975a}
\begin{align}
   \ch V_{\grpO}^\lambda &=o_\lambda=[\lambda] = \{\lambda/C\}
   =s_{\lambda/C}=\sum_{\gamma\in\C}(-1)^{|\gamma|/2} s_{\lambda/\gamma}\,,\\
   \ch V_{\grpSp}^\lambda &=sp_\lambda=\la\lambda\ra = \{\lambda/A\}
   =s_{\lambda/A}=\sum_{\alpha\in\A}(-1)^{|\alpha|/2} s_{\lambda/\alpha}\,.
\label{Eq-univOSp}
\end{align} 
These relations are the inverse of the branching rules for the restriction from
$\grpGL$ to its subgroups $\grpO$ and $\grpSp$:
\begin{align}
   \ch V_{\grpGL}^\lambda &=\{\lambda\} = [\lambda/D]
   =\sum_{\delta\in\D}\ [\lambda/\delta]
   =\sum_{\delta\in\D,\zeta\in\F}\ c^\lambda_{\delta,\zeta}\ \ch V_{\grpO}^\zeta\,,\\
   \ch V_{\grpGL}^\lambda &=\{\lambda\} = \la\lambda/B\ra
   =\sum_{\beta'\in\D}\ \la\lambda/\beta\ra 
   =\sum_{\beta'\in\D,\zeta\in\F}\ c^\lambda_{\beta,\zeta}\ \ch V_{\grpSp}^\zeta\,.
\label{Eq-GLtoOSp}
\end{align} 
That the above pairs of relations are mutually inverse is a simple consequence of the 
identities $AB=CD=1$.

In describing the Hopf algebras of the character rings of the groups
$\grpGL$, $\grpO$ and $\grpSp$ we deal only with the universal characters,
their restriction to the finite $N$ case necessitates the use of modification
rules if the relevant partitions are of too great a length. Further
details may be found elsewhere, 
for example~\cite{newell:1951a,king:1971a,black:king:wybourne:1983a,koike:terada:1987a}.

The Hopf algebra of symmetric functions, $\SY\Lambda$, is the universal, graded, bicommutative, 
biassociative self-dual Hopf algebra. 
Its properties have been spelt out in the Schur function basis in Section~\ref{Subsec-Hopf}.
Having identified in Section~\ref{Subsec-univ} the universal characters of the classical
groups and expressed them in terms of Schur functions, the Hopf algebras of
their universal character rings may be found as isomorphic copies of
$\SY\Lambda$. Despite the fact that the structure maps acting on the character
ring Hopf algebra $\CR\grpGL$, $\CR\grpO$ and $\CR\grpSp$ are isomorphic to
those of $\SY\Lambda$, they take different explicit forms in the different
canonical bases. These bases are distinguished by the use of different Littlewood
parentheses, $\{\lambda\}$, $[\lambda]$ and $\la\lambda\ra$ together with
their particular Schur-Hall scalar products with respect to which the 
bases are orthogonal.

\subsection{The general linear case}
\label{Subsec-gl}

By virtue of the identification (\ref{Eq-charGL}), the Hopf algebra, $\CR\grpGL$, of the universal character of $\grpGL$ is  immediately seen to be isomorphic to $\SY\Lambda$. Its structure is well known 
(for references see~\cite{fauser:jarvis:2003a,fauser:jarvis:king:wybourne:2005a}) and some of its
properties are summarized as follows. 

\mybenv{Theorem}\label{The-CRGL}
The ring of universal characters of $\grpGL$ is a graded self dual,
bicommutative Hopf algebra, which we denote by $\CR\grpGL$. Its 
structure maps are given by:
\begin{align}
&\textrm{product}   & 
m(\{\mu\}\otimes\{\nu\})&=\{\mu\}\cdot \{\nu\} = \{\mu\cdot\nu\}=\sum_{\zeta} 
c^\lambda_{\mu,\nu} \{\lambda\} \nn[-1ex]{}
&\textrm{unit}      &
\eta(1) &=\{0\}~~\hbox{with}~~\ \{0\}\cdot \{\mu\} 
 = \{\mu\} = \{\mu\}\cdot \{0\} \nn\nn[-1ex]{}
&\textrm{coproduct} &
\Delta(\{\lambda\}) 
&=\sum_{\mu,\nu}  
c^\lambda_{\mu,\nu} \{\mu\}\ot \{\nu\} \nn[-0.5ex]{}
&\textrm{counit}    &
\epsilon(\{\mu\}) &= \la\{0\}\mid\{\mu\} \ra =\delta_{0,\mu} \nn\nn[-1ex]{}
&\textrm{antipode}  &
S(\{\lambda\}) 
&= (-1)^{\vert\lambda\vert} \{\lambda^\prime\} \nn\nn[-1ex]{}
&\textrm{scalar product}  &
\la\cdot\mid\cdot\ra(\{\mu\}\ot\{\nu\})
&=\la\mu\mid\nu\ra=\delta_{\mu,\nu}
\end{align}
where the coefficients $c^\lambda_{\mu,\nu}$ are the Littlewood-Richardson
coefficients, $\lambda^\prime$ is the conjugate (transposed) partition and 
$\la\cdot\mid\cdot\ra : \Lambda\ot \Lambda \rightarrow\openZ$ is the usual
Schur-Hall scalar product in terms of which we have
\begin{align}
\textrm{self-duality}  &&
\la\Delta(\{\lambda\})\mid \{\mu\}\ot \{\nu\}\ra
&=\la\{\lambda\}\mid \{\mu\}\cdot\{\nu\}\ra  \,.
\end{align}
\myeenv

Because of its importance in what follows we map the antipode identity
(\ref{Eq-antipode-identity}) of $\SY\Lambda$, into the antipode identity of
$\CR\grpGL$:
\begin{align}
\sum_\nu\ (-1)^{|\nu|}\,\{\lambda/\nu\}\cdot \{\nu'\}
&=\delta_{\lambda0}\, \{0\}\,.
\label{Eq-GLai}
\end{align}

\subsection{The orthogonal case}
\label{Subsec-orth}

Having shown that the irreducible universal characters $[\lambda]$ of 
the orthogonal group $\grpO$ can be expressed in terms of universal 
characters of $\grpGL$ by $[\lambda]=\{\lambda/C\}$, 
it is possible to exploit infinite Schur function series and 
the Hopf algebra $\CR\grpGL$ to identify the action
of the structure maps on the ring of characters $[\lambda]$ forming
the canonical basis of $\CR\grpO$. This action, as will be proved 
in the following section, takes the following form:
\medskip

\mybenv{Theorem}\label{The-CRO}
The algebra $\CR\grpO$ generated by the universal characters $[\lambda]$
of the orthogonal group $\grpO$ is a bicommutative Hopf algebra. Its 
structure maps are given by:
\begin{align}
&\textrm{product}   & 
m([\mu]\cdot[\nu]) &= [\mu]\cdot[\nu] = \sum_{\zeta} [\mu/\zeta\cdot \nu/\zeta] \nn[-1ex]{}
&\textrm{unit}      &
\eta(1)&= [0]~~\hbox{with}~~\ [0]\cdot[\mu]= [\mu] = [\mu]\cdot [0] \nn\nn[-1ex]{}
&\textrm{coproduct} &
\Delta([\lambda]) 
&=\sum_{\zeta} [\lambda/(\zeta D)] \ot [\zeta]
 =\sum_{\zeta} [\lambda/\zeta]\ot [\zeta/D] 
 \nn[-1ex]{}
&\textrm{counit}    &
\epsilon([\lambda]) &= \sum_{\gamma\in C} (-1)^{\mid \gamma \mid/2} 
   \delta_{\lambda,\gamma} 
 = \delta_{\lambda,C} \nn[-1ex]{}  
&\textrm{antipode}  &
S([\lambda]) 
&= (-1)^{\vert\lambda\vert} [\lambda^\prime/(AD)] \nn \nn[-1ex]{}
&\textrm{scalar product}  &
\la\cdot\mid\cdot\ra_2([\mu]\ot[\nu])
&=\la\mu\mid\nu\ra_2=\delta_{\mu,\nu}\,.
\end{align}
\myeenv

\subsection{The symplectic case}
\label{Subsec-sp}

In the same way, by exploiting the fact that the irreducible (or
indecomposable)  universal characters $\la\lambda\ra$ of the symplectic group
$\grpSp$ can be  expressed in terms of universal characters of $\grpGL$ by 
$\la\lambda\ra=\{\lambda/A\}$, we can  identify the action of the structure maps
on the ring of characters $\la\lambda\ra$  forming the canonical basis of
$\CR\grpSp$. This action, as will be proved  in the following section, takes the
following form:

\mybenv{Theorem}\label{The-CRSp}
The algebra $\CR\grpSp$ generated by the universal characters $\la\lambda\ra$
of the symplectic group $\grpSp$ is a bicommutative Hopf algebra. Its 
structure maps are given by:
\begin{align}
&\textrm{product}   & 
m(\la\mu\ra\cdot\la\nu\ra) 
  &= \la\mu\ra\cdot\la\nu\ra = \sum_{\zeta} \la\mu/\zeta\cdot \nu/\zeta\ra \nn[-1ex]{}
&\textrm{unit}      &
\eta(1)&=\la0\ra~~\hbox{with}~~\
\la 0\ra\cdot \la\mu\ra = \la\mu\ra = \la\mu\ra\cdot \la0\ra \nn\nn[-1ex]{}
&\textrm{coproduct} &
\Delta(\la\lambda\ra) 
&=\sum_{\zeta} \la\lambda/(\zeta B)\ra \ot \la\zeta\ra
 =\sum_{\zeta} \la\lambda/\zeta\ra \ot \la\zeta/B\ra \nn[-1ex]{}
&\textrm{counit}    &
\epsilon(\la\lambda\ra) &= \sum_{\alpha\in A} (-1)^{\mid\alpha \mid/2}\, 
   \delta_{\lambda,\alpha}
  =\delta_{\lambda,A} \nn[-1ex]{}  
&\textrm{antipode}  &
S(\la\lambda\ra) 
&= (-1)^{\vert\lambda\vert} \la \lambda^\prime/(BC)\ra \nn\nn[-1ex]{}
&\textrm{scalar product}  &
\la\cdot\mid\cdot\ra_{11}(\la\mu\ra\ot\la\nu\ra)
&=\la\mu\mid\nu\ra_{11}=\delta_{\mu,\nu}\,.
\end{align}
\myeenv

\subsection{Directory of results}
\label{Subsec-dir}

All the above results, and a considerable amount of additional information,
regarding the three Hopf algebras of character rings, $\CR\grpGL$, $\CR\grpO$
and $\CR\grpSp$, are gathered together in Table~\ref{tab1}.

\begin{table}
\caption{Hopf Algebras of Character Rings, Bases and Morphisms}
\label{tab1}
\begin{align} 
\rotateleft{
\begin{array}{c|ccc}
\Lambda & \CR\grpGL & \CR\grpO & \CR\grpSp \\[1ex]
\hline
\\[1ex]
1&\{0\}&[0]&\la0\ra
\\[1ex]
p_n & \sum_{b=0}^{n-1} (-1)^b\{n-b,1^b\} 
    & \sum_{b=0}^{n-1} (-1)^b[n-b,1^b]+\chi(2\vert n)[0]
    & \sum_{b=0}^{n-1} (-1)^b\la n-b,1^b\ra +\chi(2\vert n)\la 0\ra 
    \\[1ex] 
h_n & \{n\}
    & [n/D]
    & \la n\ra 
    \\[1ex]
e_n & \{1^n\}
    & [1^n]
    & \la 1^n/B\ra 
    \\[1ex]
s_\lambda & \{\lambda\}
    & [\lambda/D]
    & \la \lambda/B\ra 
    \\[1ex]
o_\lambda &\{\lambda/C\}
    & [\lambda]
    & \la \lambda/(BC)\ra 
    \\[1ex]
sp_\lambda &\{\lambda/A\}
    & [\lambda/(AD)]
    & \la \lambda\ra 
    \\[1ex]
\hline 
    \\[1ex]     
m   & m(\{\mu\}\ot \{\nu\})=\{\mu\cdot\nu\}
    & m([\mu]\ot [\nu]) = \sum_{\zeta\in\P} [\mu/\zeta\cdot \nu/\zeta]
    & m(\la\mu\ra\ot\la\nu\ra) = \sum_{\zeta\in\P}\la\mu/\zeta\cdot \nu/\zeta\ra
    \\[1ex]
\Delta & \Delta(\{\lambda\})=\sum_{\zeta\in\P} \{\lambda/\zeta\}\ot\{\zeta\}
    &\Delta([\lambda])=\sum_{\zeta\in\P} [\lambda/\zeta]\ot [\zeta/D]
    &\Delta(\la\lambda\ra)=\sum_{\zeta\in\P} \la\lambda/\zeta\ra\ot \la\zeta/B\ra
    \\[1ex]    
\epsilon & \epsilon(\{\lambda\})=\delta_{\lambda,0}
    & \epsilon([\lambda])=\sum_{\gamma\in C}(-1)^{|\gamma|/2}\delta_{\lambda,\gamma}
    & \epsilon(\la\lambda\ra)=\sum_{\alpha\in A}(-1)^{|\alpha|/2}\delta_{\lambda,\alpha}
    \\[1ex]
\eta & \eta(1)=\{0\}
    & \eta(1)=[0]
    & \eta(1)=\la 0\ra
    \\[1ex]
S   & S(\{\lambda\}) = (-1)^{\vert\lambda\vert}\{\lambda^\prime\}
    & S([\lambda]) = (-1)^{\vert\lambda\vert} [\lambda^\prime/(AD)]
    & S(\la\lambda\ra) 
       = (-1)^{\vert\lambda\vert}\la\lambda^\prime/(BC)\ra
    \\[1ex]
\hline 
\\[1ex]
\la\cdot|\cdot\ra&\la\{\lambda\}\mid\{\mu\}\ra=\delta_{\lambda,\mu}
                 &\la[\lambda]\mid[\mu]\ra_{2}=\delta_{\lambda,\mu}
                 &\la\la\lambda\ra\mid\la\mu\ra\ra_{11}=\delta_{\lambda,\mu}
                 \\[1ex]
\hline                 
\end{array}
} 
\nonumber
\end{align}
\end{table}

The first column of this directory gives the abstract Hopf algebra notation  for
bases and morphisms of $\SY\Lambda$, for any $n\in\openN$ and $\lambda\in\P$. 
The second column gives the notion for the
Hopf algebra of the universal character ring of the general linear  group, as
studied for example in~\cite{fauser:jarvis:2003a}. The third and fourth columns
provide the isomorphic images of the structure maps and bases in the character
rings of the orthogonal and symplectic groups. We have used the notational convention
whereby $\chi(P)$ is the truth symbol, that is $\chi(P)=1$ if the proposition $P$ is true, 
and $0$ otherwise. Thus $\chi(2\vert n)=1$ if $n$ is even and $\chi(2\vert n)=0$ if $n$ 
is odd.

\paragraph{Remark.} 
While $\Lambda$ and $\CR\grpGL$ share the same Schur-Hall scalar
product we emphasise that $\CR\grpO$ and $\CR\grpSp$ can quite naturally be equipped 
with \emph{new structure maps}, plethystic Schur-Hall scalar products, indexed by $2$ and
$11$, which are defined as shown in Table~\ref{tab1} so as to ensure that the orthogonal and 
symplectic Schur functions form orthonormal bases of $\CR\grpO$ and $\CR\grpSp$, respectively.\qed 
\medskip

The precise definitions of the bases involved in some of the formulae of Table~\ref{tab1}
will be given in the following sections. However, this table makes it clear that there
are \emph{unique} instances of symmetric functions, such as power sum symmetric
functions, which are tied to the underlying alphabet and are, up to isomorphism,
equivalent in all the character Hopf algebras under consideration. Despite this, if
written in the canonical basis of a specific character Hopf algebra, it can be
seen that such objects may look different and may also exhibit combinatorial
differences.

\section{Orthogonal and symplectic character ring Hopf algebras}
\label{Sec-Hopf-proofs}

In this section we provide proofs of the validity of each of the structure map 
formulae listed in Table~\ref{tab1}.

\subsection{The case of \CR$\grpO$}
\label{Subsec-o-Hopf}

We consider in turn each of the structure maps listed in Theorem~\ref{The-CRO}.

The {\bf product} formula
\begin{align}
[\mu]\cdot[\nu] = \sum_{\zeta\in\P} [\mu/\zeta\cdot \nu/\zeta]
\label{Eq-NL-O}
\end{align}
is a classical result of Newell~\cite{newell:1951a} and Littlewood~\cite{littlewood:1958a} 
that appears as a special case of the development 
in~\cite{fauser:jarvis:king:wybourne:2005a} for more general subgroups of 
the general linear group. Its derivation can be accomplished most easily by noting
that $[\mu]\cdot[\nu]=\{\mu/C\}\cdot\{\nu/C\}=[(\{\mu/C\}\cdot\{\nu/C\})/D]$ 
where the coefficient of $\{\lambda\}$ in $(\{\mu/C\}\cdot\{\nu/C\})/D$ is given by
\begin{align}
\la \{\lambda\} \,|\, (\{\mu/C\}\cdot\{\nu/C\})/D \ra
&= \la \{\lambda\}\cdot D \,|\, \{\mu/C\}\cdot \{\nu/C\} \ra
= \la \Delta(\{\lambda\}\cdot D)\,|\, \{\mu/C\}\otimes \{\nu/C\} \ra\cr\cr
&= \sum_{\zeta\in\P} \la \Delta(\{\lambda\}) \cdot (D\otimes D) \cdot \{\zeta\}\otimes\{\zeta\}\,|\, \{\mu/C\}\otimes \{\nu/C\} \ra\cr
&= \sum_{\zeta\in\P} \la \Delta(\{\lambda\})\,|\, \{\mu/(CD\zeta)\}\otimes \{\nu/(CD\zeta)\} \ra\cr
&= \sum_{\zeta\in\P} \la \Delta(\{\lambda\})\,|\, \{\mu/\zeta\}\otimes \{\nu/\zeta\} \ra
= \sum_{\zeta\in\P} \la \{\lambda\}\,|\, \{\mu/\zeta\}\cdot\{\nu/\zeta\} \ra\,,
\end{align}
from which the result (\ref{Eq-NL-O}) follows.
 
To find the {\bf coproduct} we need to find first the ordinary coproduct of
a skew Schur function. This can be looked up in Macdonald~\cite{macdonald:1979a} (Eq. 5.9 and 5.10, p72). 
The idea is to expand 
$s_\lambda(x,y,z)$, a double coproduct, in two different ways:
\begin{align}
s_\lambda(x,y,z)
&= \sum_\nu s_{\lambda/\nu}(x,y)s_\nu(z) \nn
&= \sum_\mu s_{\lambda/\mu}(x)s_\mu(y,z) 
 = \sum_{\mu,\nu} s_{\lambda/\mu}(x)s_{\mu/\nu}(y)s_\nu(z)\,,
\end{align}
and comparing coefficients of $s_\nu(z)$ gives  
\begin{align}
s_{\lambda/\nu}(x,y)
&=\sum_{\mu} s_{\lambda/\mu}(x)s_{\mu/\nu}(y)\,,
&&~~\hbox{that is}~~&
\Delta(s_{\lambda/\nu})
&=\sum_{\mu} s_{\lambda/\mu}\ot s_{\mu/\nu}\,.
\end{align}
Now we can proceed to compute
\begin{align}
\Delta([\lambda]):&= \Delta\{\lambda/C\}
=\sum_{\gamma\in C} (-1)^{|\gamma|/2}\Delta(\{\lambda/\gamma\})
=\sum_{\gamma\in C,\zeta\in\P} (-1)^{|\gamma|/2}\{\lambda/\zeta\}\ot \{\zeta/\gamma\} \nn
&=\sum_{\zeta\in\P} \{\lambda/\zeta\}\ot [\zeta]
=\sum_{\zeta\in\P} [(\lambda/\zeta)/D]\ot [\zeta] 
=\sum_{\zeta\in\P} [\lambda/(\zeta D)] \ot [\zeta]\,.
\label{Eq-DeltaO-1}
\end{align}
This can equally well be rewritten to give a second form of the coproduct 
derived using a pair of related expansions of a skew Schur function,
$s_{\lambda/\zeta}=\sum_\sigma c^\lambda_{\zeta,\sigma} s_\sigma$ and
$s_{\lambda/\sigma}=\sum_\zeta c^\lambda_{\zeta,\sigma} s_\zeta$, a move we use
below frequently. Here it gives,
\begin{align}
\Delta([\lambda])
&=\sum_{\zeta\in\P} [\lambda/(\zeta D)] \ot [\zeta]=\sum_{\zeta\in\P} [(\lambda/\zeta)/D)] \ot [\zeta] \nn
&=\sum_{\zeta,\sigma\in\P} c^\lambda_{\zeta,\sigma}\ [\sigma/D]\ot [\zeta]=\sum_{\sigma\in\P} [\sigma/D]\ot [\lambda/\sigma]\,.
\label{Eq-DeltaO-2}
\end{align}
The coproduct $\Delta([\lambda])$ is cocommutative, as can be seen by using in
the same way as above the connection between the outer product of Schur
functions $s_\zeta\cdot s_\delta=\sum_\eta c^\sigma_{\zeta,\delta}\,s_\sigma$
and the skew Schur functions expansion $s_{\sigma/\delta}=\sum_\zeta
c^\sigma_{\zeta,\delta}\,s_\zeta$, to obtain a third form:
\begin{align}
\Delta([\lambda])
&=\sum_{\zeta\in\P} [\lambda/(\zeta D)] \ot [\zeta]=\sum_{\zeta\in\P,\delta\in D} [(\lambda/(\zeta\cdot\delta)] \ot [\zeta]\nn
&=\sum_{\zeta,\sigma\in\P,\delta\in D} c^\sigma_{\zeta,\delta}\ [\lambda/\sigma] \ot [\zeta]
=\sum_{\sigma\in\P,\delta\in D} [\lambda/\sigma]\ot [\sigma/\delta]=\sum_{\sigma\in\P} [\lambda/\sigma]\ot [\sigma/D]\,.
\label{Eq-DeltaO-3}
\end{align}
Finally, we may use $s_{\lambda/\sigma}=\sum_\zeta c^\lambda_{\sigma,\zeta} s_\zeta$ and 
$\sum_\sigma c^\lambda_{\sigma,\zeta} s_\sigma= s_{\lambda/\zeta}$ to obtain the
fourth form
\begin{align}
\Delta([\lambda])
&=\sum_{\sigma\in\P} [\lambda/\sigma]\ot [\sigma/D]=\sum_{\sigma,\zeta\in\P} c^\lambda_{\sigma,\zeta} [\zeta] \ot [\sigma/D]\nn
&=\sum_{\zeta\in\P} [\zeta] \ot [(\lambda/\zeta)/D]
=\sum_{\zeta\in\P} [\zeta]\ot [\lambda/(\zeta D)]\,.
\label{Eq-DeltaO-4}
\end{align}

The actions of the {\bf counit}, $\epsilon$, the {\bf unit}, $\eta$, and 
the {\bf antipode}, $S$, follow immediately from 
their action in \CR$\grpGL$ and the fact that $[\lambda]=\{\lambda/C\}$.
Thus in the \CR$\grpO$ basis  
\begin{align}
  \epsilon([\lambda])&=\epsilon(\{\lambda/C\})=\sum_{\gamma\in\C} (-1)^{|\gamma|/2} \epsilon(\{\lambda/\gamma\})
  =\sum_{\gamma\in\C} (-1)^{|\gamma|/2} \delta_{\lambda,\gamma}\,;\cr
  \eta(1)&=\{0\}=[0/D]=[0]\,; \cr\cr 
  S([\lambda])&=S(\{\lambda/C\}=\sum_{\gamma\in\C} (-1)^{|\gamma|/2} S(\{\lambda/\gamma\})
   =\sum_{\gamma\in\C} (-1)^{|\gamma|/2} (-1)^{|\lambda|-|\gamma|} \{\lambda'/\gamma'\}\cr
  &=(-1)^{|\lambda|} \{\lambda/C'\} =(-1)^{|\lambda|} \{\lambda/A\} = (-1)^{|\lambda|} [(\lambda/A)/D]
   = (-1)^{|\lambda|} [\lambda/(AD)]
\end{align} 
all as shown in Table~~1.

By exploiting the Schur-Hall scalar product we may reinterprit $\epsilon$ and introduce
its convolutive inverse $\epsilon^{-1}$ as follows: 
\mybenv{Definition}
The counit $\epsilon$ and its convolutive inverse $\epsilon^{-1}$ for \CR$\grpO$ 
may be interpreted as linear forms $c$ and $d$: \CR$\grpGL\rightarrow\openZ$ 
defined as follows:%
\footnote{This definition should be compared with a slightly different point of
view developed in the section on adapted normal ordered products 
in~\cite{brouder:fauser:frabetti:oeckl:2002a}, which can be used to define a
quantum field theory on an external background.}
\begin{align}
&\epsilon([\lambda])=c(\{\lambda\})
~~\hbox{with}~~c(\{\lambda\}):=\la C \mid \{\lambda\}\ra
=\sum_{\gamma\in\C} (-1)^{|\gamma|/2} \la \{\gamma\} \mid \{\lambda\}\ra
=\sum_{\gamma\in\C} (-1)^{|\gamma|/2} \delta_{\gamma,\lambda} \,,\nn
\\
&\epsilon^{-1}([\lambda])=d(\{\lambda\})
~~\hbox{with}~~d(\{\lambda\}):=\la D \mid \{\lambda\}\ra
=\sum_{\delta\in\D} \la \{\delta\} \mid \{\lambda\}\ra 
=\sum_{\delta\in\D} \delta_{\delta,\lambda} \,.\nonumber
\end{align}
\label{Def-cd}
\myeenv
\mybenv{Corollary} (see~\cite{fauser:jarvis:2003a})
The linear forms (1-cochains)  $c$ and $d$ are convolutive inverses with
respect to the $\CR\grpGL$ outer coproduct and product in $\openZ$.
\label{Cor-cd}
\myeenv
\paragraph{Proof:}
\begin{align}
(c\conv d)(\{\lambda\}) 
&=\sum_{(\lambda)} c(\{\lambda_{(1)}\})d(\{\lambda_{(2)}\}) 
=\sum_{(\lambda)} \la C\mid \{\lambda_{(1)}\}\ra 
                   \la D\mid \{\lambda_{(2)}\}\ra \nn
&=\sum_{(\lambda)} \la C\ot D \mid \{\lambda_{(1)}\}\ot \{\lambda_{(2)}\}\ra 
= \la  C\ot D \mid \Delta(\{\lambda\}) \ra\nn
&=\la CD \mid \{\lambda\}\ra = \la 1 \mid \{\lambda\}\ra 
=\epsilon(\{\lambda\}) = \delta_{\lambda,0}\,.
\end{align}
\qed

Finally, we might check that the product and coproduct are mutual 
{\bf  coalgebra} and {\bf algebra homomorphisms}. We establish this fact
by direct computation:
\begin{align}
&(\Delta\, m )([\lambda]\ot [\mu])
=\sum_{\zeta}\Delta([\lambda/\zeta\cdot \mu/\zeta]) 
=\sum_{\rho,\zeta}
  [\left( \lambda/\zeta\cdot \mu/\zeta\right)/\rho]\ot [\rho/D] \nn
&=\sum_{\sigma,\rho,\zeta}
  [\lambda/(\zeta\sigma)\cdot \mu/(\zeta(\rho/\sigma))]
  \ot [\rho/D] 
=\sum_{\xi,\sigma,\zeta}
  [\lambda/(\zeta\sigma)\cdot \mu/(\zeta\xi)]
  \ot [(\sigma\xi)/D] \nn
&=\sum_{\tau,\xi,\sigma,\zeta}
  [\lambda/(\zeta\sigma)\cdot \mu/(\zeta\xi)]
  \ot [\sigma/(\tau D)\cdot \xi/(\tau D)] 
=\sum_{\sigma,\xi}
  \left( [\lambda/\sigma]\cdot [\mu/\xi] \right)
  \ot 
  \left( [\sigma/D] \cdot [\xi/D] \right) \nn
&=\sum_{\sigma,\xi}
  \left( [\lambda/\sigma]\ot [\sigma/D] \right)
  \cdot
  \left( [\mu/\xi]\ot [\xi/D]\right) 
=\Delta([\lambda]) \cdot \Delta([\mu]) = m(\Delta([\lambda])\otimes \Delta([\mu]))\,,
\end{align}
showing the claim.\qed
\paragraph{Remarks.}
It could be argued that the above proof is unnecessary. We considered just a
linear isomorphism on the module underlying the symmetric function Hopf algebra,
and the result is in a natural way, a homomorphic image. However, the displayed
calculations show explicitly how the structure maps are written in the
orthogonal Schur function bases, how the combinatorics alters, and that
everything is set up correctly.

Note also the most remarkable fact that the structure of the Hopf algebra
\CR$\grpO$ does not distinguish between even and odd orthogonal groups. It does
not even rely on the fact that the metric tensor $g_{ij} = g_{ji}$ of Schur
symmetry type $\{2\}$, which defines the orthogonal group, is invertible. Such
degenerate cases are instances of Cayley-Klein groups (see conclusions for
further comments). The even or oddness of the underlying group will show up
in a subtle way when we define particular bases for these Hopf algebras below.

\subsection{The case of $\CR\grpSp$}
\label{Subsec-sp-Hopf}

The validity of the structure maps of $\CR\grpSp$ given in
Theorem~\ref{The-CRSp}  may be established by copying and pasting the proof for
the orthogonal case.  One merely changes all orthogonal characters into
symplectic ones,   $[\lambda]\rightarrow\la\lambda\ra$, and interchanges Schur
function series, $C\leftrightarrow A$ and $D\leftrightarrow B$. All arguments run
through as before. In the case of the counit, it is also necessary to
interchange the labelling on the linear forms, $c\rightarrow a$ and 
$d\rightarrow b$, where by analogy with Definition~\ref{Def-cd} we have:

\mybenv{Definition}
The counit $\epsilon$ and its convolutive inverse $\epsilon^{-1}$ for \CR$\grpSp$ 
may be interpreted as linear forms $a$ and $b$: \CR$\grpGL\rightarrow\openZ$ 
defined as follows:%
\begin{align}
&\epsilon(\la\lambda\ra)=a(\{\lambda\})
~~\hbox{with}~~a(\{\lambda\}):=\la A \mid \{\lambda\}\ra
=\sum_{\alpha\in A} (-1)^{|\alpha|/2} \la \{\alpha\} \mid \{\lambda\}\ra
=\sum_{\alpha\in A} (-1)^{|\alpha|/2} \delta_{\alpha,\lambda} \,,\nn
\\
&\epsilon^{-1}(\la\lambda\ra)=b(\{\lambda\})
~~\hbox{with}~~b(\{\lambda\}):=\la B \mid \{\lambda\}\ra
=\sum_{\beta\in B} \la \{\beta\} \mid \{\lambda\}\ra 
=\sum_{\beta\in B} \delta_{\beta,\lambda} \,.\nonumber
\end{align}
\label{Def-ab}
\myeenv

Once again as in Corollary~\ref{Cor-cd} we have: 
\mybenv{Corollary} (see~\cite{fauser:jarvis:2003a})
The linear forms (1-cochains) $a$ and $b$ are convolutive inverses with
respect to the $\CR\grpGL$ outer coproduct and product in $\openZ$.
\label{Cor-ab}
\myeenv

\section{Bases for $\CR\grpO$ and $\CR\grpSp$}
\label{Sec-bases}

\subsection{Power sum symmetric functions}
\label{Subsec-powsum}

A major issue in setting the above abstract machinery to work in concrete
(physical) examples, is a proper identification in the various character rings
of the usual canonical bases of the symmetric function ring.  In making this
identification, we will encounter some familiar and also some surprising
results. We start with the power sum symmetric functions on a  finite number of
variables $N$. The one part power sum symmetric functions are  defined on the
variables $(x_1,\ldots,x_N)$ by
\begin{align}
p_n &:=\sum_{i=1}^N x_i^n \,
\end{align}
which is independent of the meaning of the alphabet. 

In the $\grpGL(N)$ case
the $x_i$ are the eigenvalues of a $\grpGL(N)$ element $g$ within a
$\grpGL(N)$ conjugacy class. There is the constraint $\prod_i x_i \ne 0$ in
force to ensure invertibility. We use the well known hook expansion in terms of 
the Schur functions identified with irreducible $\grpGL(N)$ characters:
\begin{align}
p_n(x_1,\ldots,x_n) & =\sum_{a+b+1=n} (-1)^b \{a+1,1^b\}(x_1,\ldots,x_N)
\end{align}
This formula is stable with respect to the limit $N\rightarrow\infty$
so that we immediately have in the case of $\CR\grpGL$ the identification
\begin{align}
p_n & =\sum_{a+b+1=n} (-1)^b \{a+1,1^b\} \,.
\label{Eq-pnGL}
\end{align}

In branching to orthogonal or the symplectic groups, as one can see from
(\ref{Eq-chVuniv}) the eigenvalues now generally speaking come in pairs $x_k$ and $\xb_k$
and we can split $p_n$ into at least two parts. 
In the orthogonal $O(N)$ case, there are four possibilities, and in the symplectic
$Sp(N)$ case there are two. Confining attention to the unimodular case it follows
from (\ref{Eq-chVuniv}) that:
\begin{align}
p_n(x,\xb)
&=p_n(x)+p_n(\xb)   &&\hbox{for}~~\grpSO(2K)\,;\nn
p_n(x,\xb,1) 
&=p_n(x)+p_n(\xb)+1 &&\hbox{for}~~\grpSO(2K+1)\,;\nn
p_n(x,\xb) 
&=p_n(x)+p_n(\xb)   &&\hbox{for}~~\grpSp(2K)\,;\nn
p_n(x,\xb,1) 
&=p_n(x)+p_n(\xb)+1 &&\hbox{for}~~\grpSp(2K+1)\,,
\label{Eq-pn}
\end{align}
where in each case $p_n(x)=\sum_{i=1}^K x^n_i$ and 
$p_n(\xb)=\sum_{i=1}^K \xb^n_i$.

Clearly, this is the place where the dimensionality $N=2K$ or $N=2K+1$ comes
into play. However, this does not prevent us from establishing a result stable
in the $K\rightarrow\infty$ limit. Indeed we find as a corollary to
(\ref{Eq-pnGL}) the result  appropriate to $\CR\grpO$:
\mybenv{Corollary}
\begin{align}
p_0=[0]~~~~\hbox{and}~~~~p_n &=\sum_{a+b+1=n} (-1)^b [a+1,1^b] + \chi(2\vert n)[0]
&&\hbox{for}~~n\ge1\,,
\end{align}
where $\chi$ is the truth function, so that $\chi(2\vert n)=1$ if $n$ is
even and $\chi(2\vert n)=0$ if $n$ is odd.  
\myeenv

Note that the $n$ in this truth function has to do
with the index of the one part power sums, and not with the number of
its variables! 

\paragraph{Proof:}
For $n=0,1$ we can directly verify that $p_0=[0]$ and $p_1=[1]$, thereby
proving the statement in these cases. Henceforth we assume $n\ge2$.
The $D$ series partitions $\delta\in\D=2\P$ have only even parts, 
and of these only the partitions of type $\{2k\}$ can fit into a hook. Thus
\begin{align}
p_n
&=\sum_{a+b+1=n} (-1)^{b} \{a+1,1^b\}
 =\sum_{b\ge0}^{n-1} (-1)^{b} \{n-b,1^b\}\nn
&=\sum_{b\ge0}^{n-1} (-1)^{b} [(n-b,1^b)/D]
 =\sum_{k\ge0}^{[n/2]}\sum_{b\ge0}^{n-1} (-1)^{b} [(n-b,1^b)/(2k)] \nn
&=\sum_{b\ge0}^{n-1} (-1)^{b} [(n-b,1^b)]
 +\sum_{k\ge1}^{[n/2]}\sum_{b\ge0}^{n-2k} (-1)^{b} [(n-b-2k)\cdot(1^b)] \nn
&=\sum_{b\ge0}^{n-1} (-1)^{b} [(n-b,1^b)]
 +\sum_{k\ge1}^{[n/2]}\sum_{\zeta\in\P} (-1)^{\vert\zeta\vert} 
   [(n-2k)/\zeta \cdot \zeta' ] \nn
&=\sum_{b\ge0}^{n-1} (-1)^{b} [(n-b,1^b)] + \sum_{k\ge1}^{[n/2]}\delta_{n-2k,0}[0] \nn
&=\sum_{b\ge0}^{n-1} (-1)^{b} [(n-b,1^b)] + \chi(2\vert n)[0] \,,
\label{Eq-pn-O}
\end{align}
where, in the penultimate line, the second term has resulted
from the antipode property (\ref{Eq-GLai}).
\qed

We have an exactly analogous result for $\CR\grpSp$:
\mybenv{Corollary}
\begin{align}
p_0=\la0\ra~~~~\hbox{and}~~~~p_n &=\sum_{a+b+1=n} (-1)^b \la a+1,1^b\ra + \chi(2\vert n)\la0\ra
&&\hbox{for}~~n\ge1\,.
\end{align}
\vskip-1ex
\myeenv

We recall that the one part power sums are the primitive elements of the
symmetric function Hopf algebra. They form a rational basis of this Hopf
algebra. This implies:
\mybenv{Proposition}
The one part power sums $p_n$ map to the primitive elements of the Hopf algebra
of the universal character rings of $\grpGL$, $\grpO$ and $\grpSp$.  That is, in
each case we have
\begin{align}
\Delta(p_n) &= p_n\ot 1 + 1\ot p_n \,. 
\label{Eq-Delta-pn}
\end{align}
\vskip-2ex
\myeenv
\vskip-2ex

\paragraph{Proof:}
This is a trivial consequence of (\ref{Eq-coprod-hep}), since the isomorphism of
Hopf algebras which we have established is  independent of the underlying
alphabet, and hence  does not alter the coproduct properties of the power sums.
\qed

\subsection{Complete symmetric functions}
\label{Subsec-complete}

In this section we investigate the nature of complete symmetric functions in
each of our three rings of universal characters by means of the
maps between Schur functions and the characters.

First, our maps allow us to see immediately that, in accordance with
the formulae of Table~1, we have
\begin{align}
 h_n&=s_n=\{n\}\,,\nn
 h_n&=s_n=\{n\}=[n/D]=\sum_k\ [n/(2k)]=\sum_{k=0}^{[n/2]} [n-2k]\,,\nn
 h_n&=s_n=\{n\}=\la n/B\ra= \la n\ra,
\label{Eq-hnCR} 
\end{align}
where $[n/2]$ is the integer part of $n/2$. Moreover, we have
\mybenv{Proposition}
The above images of the one part complete symmetric functions $h_n$ under the
maps from the Hopf algebra of symmetric functions to the universal character 
rings of $\grpGL$, $\grpO$ and $\grpSp$ are divided 
powers~\cite{newman:1972a,berkson:newman:1978a,sweedler:1979a}, their
coproducts take the form:
\begin{align}
 \Delta(h_n) &= \sum_r\ h_{n-r}\ot h_r.
\end{align}
\vskip-2ex
\myeenv
\vskip-2ex

\paragraph{Proof:} These results are a direct consequence of
(\ref{Eq-coprod-hep}), since the maps between the Hopf algebras are
isomorphisms, but they can also be derived as follows.
\begin{align}
\Delta(\{n\}) &= \sum_{\zeta} \{n/\zeta\}\ot \{\zeta\} 
= \sum_r \{n/r\}\ot \{r\} = \sum_r \{n-r\}\ot\{r\}\,,\nn
\Delta([n/D]) &= \sum_{\zeta} [n/(\zeta D)]\ot [\zeta/D] 
= \sum_r [(n/r)/D]\ot [r/D] = \sum_r [(n-r)/D] \ot [r/D]\,,\nn
\Delta(\la n\ra) &= \sum_{\zeta} \la n/\zeta\ra \ot \la \zeta/B\ra 
= \sum_r \la n/r\ra \ot \la r/B\ra = \sum_r \la n-r\ra \ot \la r\ra\,.
\end{align}

\subsection{Elementary symmetric functions}
\label{Subsec-elem}

The elementary symmetric functions $e_n$ map as follows to the three character 
rings of interest:
\begin{align}
e_n &=s_{1^n}=\{1^n\}\,,\nn
e_n &=s_{1^n}=\{1^n\}=[1^n/D]=[1^n]\,, \nn
e_n &=s_{1^n}=\{1^n\}=\la 1^n/B\ra =\sum_r \la 1^{n-2r}\ra = \sum_{r=0}^{[n/2]} \la 1^{n-2r}\ra\,. 
\end{align}
Moreover, we have
\mybenv{Proposition}
The above images of the one part elementary symmetric functions $e_n$ under the
maps from the Hopf algebra of symmetric functions to the universal character 
rings of $\grpGL$, $\grpO$ and $\grpSp$ are again divided 
powers since their coproducts all take the form:
\begin{align}
 \Delta(e_n) &= \sum_r\ e_{n-r}\ot e_r.
\end{align}
\vskip-2ex
\myeenv
\vskip-2ex

\paragraph{Proof:} These results are a direct consequence of
(\ref{Eq-coprod-hep}), since the maps between the Hopf algebras are
isomorphisms, but they can also be derived as follows.
\begin{align}
\Delta(\{1^n\}) &= \sum_{\zeta} \{1^n/\zeta\}\ot \{\zeta\} 
= \sum_r \{1^n/1^r\}\ot \{1^r\} = \sum_r \{1^{n-r}\}\ot\{1^r\}\,,\nn
\Delta([1^n]) &= \sum_{\zeta} [1^n/\zeta]\ot [\zeta/D] 
= \sum_r [1^n/1^r]\ot [1^r/D] = \sum_r [1^{n-r}] \ot [1^r]\,,\\
\Delta(\la 1^n/B\ra) &=\! \sum_{\zeta} \la 1^n/(\zeta B)\ra \ot \la \zeta/B\ra 
=\! \sum_r \la (1^n/1^r)/B\ra \ot \la 1^r/B\ra 
\nn& ~~~~~~~~~~~~~~~~~~~~~~~~~~~~~~~~~~~~~~~~~~~~~~ 
=\! \sum_r \la 1^{n-r}/B\ra \ot \la 1^r/B\ra.\nonumber
\end{align}

\section{Scalar products, adjoints, Foulkes derivatives and duals}
\label{Sec-scalar-dual}

\subsection{Scalar products}
\label{Subsec-scalar}

Unlike the Schur functions of general linear type, 
it can be readily checked, that Schur
functions of orthogonal type $o_\lambda=[\lambda]$ and of the symplectic
type $sp_\lambda=\la\lambda\ra$ are {\bf not} orthogonal with respect to the
Schur-Hall scalar product. It is hence necessary to 
define new `orthogonal' and `symplectic'
scalar products, accounting for the fact that we consider the 
orthogonal and symplectic Schur functions to be universal characters of 
irreducible orthogonal and symplectic group representations. 
The orthogonality of the universal characters of \CR$\grpGL$, \CR$\grpO$
and $\CR\grpSp$ are expressed through the following: 

\mybenv{Definition}
The general linear, orthogonal and symplectic scalar products are defined by:
\begin{align}
\la \cdot \mid \cdot \ra 
: \textrm{\CR}\grpGL \ot \textrm{\CR}\grpGL\rightarrow \openZ &
~~~~\hbox{with}~~~~
\la \{\lambda\} \mid \{\mu\} \ra = \delta_{\lambda,\mu}\,;\nn
\la \cdot \mid \cdot \ra_{2} 
: \textrm{\CR}\grpO \ot \textrm{\CR}\grpO\rightarrow \openZ &
~~~~\hbox{with}~~~~
\la [\lambda] \mid [\mu]\ra_{2} = \delta_{\lambda,\mu}\,;\nn
\la \cdot \mid \cdot \ra_{11} 
: \textrm{\CR}\grpO \ot \textrm{\CR}\grpO\rightarrow \openZ &
~~~~\hbox{with}~~~~
\la \la\lambda\ra \mid \la\mu\ra\ra_{11} = \delta_{\lambda,\mu}\,,
\label{Eq-scalprod-O}
\end{align}
for all partitions $\lambda$ and $\mu$
\myeenv

\noindent 
The indices $2$ and $11$ are a reminder of the plethystic character of the
branching from $\grpGL$ to $\grpO$ (see~\cite{fauser:jarvis:king:wybourne:2005a} 
and the previous introductory remarks). 

The relation between
the scalar products of $\CR\grpO$ and $\CR\grpSp$ and those of $\CR\grpGL$ is
such that
\begin{align}
\la \[\lambda\] \mid \[\mu\] \ra_2
&=\delta_{\lambda,\mu}
= \la \{\lambda\} \mid \{\mu\} \ra 
= \la \[\lambda/D\] \mid \[\mu/D\] \ra \,;\nn
\la \la\lambda\ra \mid \la\mu\ra \ra_{11}
&=\delta_{\lambda,\mu}
= \la \{\lambda\} \mid \{\mu\} \ra 
= \la \la\lambda/B\ra \mid \la\mu/B\ra \ra \,.\nn
\end{align}

We now consider two maps from the ring of symmetric functions
$\Lambda$ into the ring $\End(\Lambda)$ and their general linear, orthogonal and symplectic
counter parts. These are the operators, $\cdot$ and $\perp$:~~$\Lambda \rightarrow \End(\Lambda)$
corresponding to `multiplying by a Schur function' and its adjoint `skewing with a Schur function', 
which we have used frequently above.
\begin{align}
&s_\lambda\cdot(s_\mu)
  =\sum_\nu c_{\lambda,\mu}^\nu s_\nu= s_{\lambda\cdot\mu} 
  \qquad\hbox{and}\qquad
 s_\lambda^\perp(s_\mu) 
  =\sum_\nu c^\mu_{\lambda,\nu} s_\nu =s_{\mu/\lambda} \,.
\end{align}
These two operations are related via the Schur-Hall 
scalar product
\begin{align}
\la s_\mu\cdot(s_\nu) \mid s_\lambda \ra 
&= \la s_\mu\cdot s_\nu \mid s_\lambda\ra
=c_{\mu,\nu}^\lambda = \la s_\nu \mid s_{\lambda/\mu}\ra  
 = \la s_\nu \mid s_\mu^\perp(s_\lambda) \ra\,.
\label{Eq-sfn-mult-perp}
\end{align}

Formulae analogous to this exist for all our classical group universal
characters:
\mybenv{Corollary}
For all $\lambda,\mu,\nu\in\P$
\begin{align}
\la \{\mu\cdot\nu\} \mid \{\lambda\}\ra &= c_{\mu,\nu}^\lambda = \la \{\nu\} \mid \{\lambda/\mu\} \ra\,; \nn
  \la [\mu\cdot\nu] \mid [\lambda] \ra_2 &= c_{\mu,\nu}^\lambda = \la [\nu] \mid [\lambda/\mu] \ra_2\,; \nn
  \la \la\mu\cdot\nu\ra \mid \la\lambda\ra \ra_{11} &= c_{\mu,\nu}^\lambda = \la \la\nu\ra \mid \la\lambda/\mu\ra \ra_{11}\,.
\end{align}
\myeenv

\paragraph{Proof:}
The first of these follows from (\ref{Eq-sfn-mult-perp}) through the usual identification
$\{\lambda\}=s_\lambda$ for all $\lambda$. For the second, one merely notes that from (\ref{Eq-scalprod-O})  
\begin{align}
  \la [\lambda] \mid [\mu\cdot\nu] \ra_2 &= \sum_\zeta\ c_{\mu,\nu}^\zeta  \la [\lambda]\mid[\zeta]\ra_2 = c_{\mu,\nu}^\lambda\,;\nn
  \la [\lambda/\mu] \mid [\nu] \ra_2 &= \sum_\zeta\ c_{\mu,\zeta}^\lambda  \la [\zeta]\mid[\nu]\ra_2 = c_{\mu,\nu}^\lambda\,.
\end{align}
The third is derived in an analogous manner.
\qed

However it should be noted that these relations do not help us identify an adjoint of multiplication
for $\CR\grpO$ and $\CR\grpSp$ since $[\mu]\cdot[\nu]\neq[\mu\cdot\nu]$ and
$\la\mu\ra\cdot\la\nu\ra\neq\la\mu\cdot\nu\ra$.

The adjoint of multiplication by a Schur function with respect to the
Schur-Hall scalar product, that is the skew or $\perp$, is called the Foulkes
derivative. This can be used to introduce differential operators, for example 
in Macdonald~\cite{macdonald:1979a} one finds both
\begin{align}\label{Eq-pnperp}
p_n^\perp &= n \frac{\partial}{\partial p_n} 
&&\hbox{and}&&
p_n^\perp=\sum_{r\ge0} h_r \frac{\partial}{\partial h_{n+r}}\,.
\end{align}
This leads to the interesting fact, that the coproduct can be written in
terms of the adjoint:
\begin{align}
\Delta(f) 
&=\sum_{\mu} s_\mu^\perp(f) \ot s_\mu 
 =\sum_{\mu,(f)} \epsilon(s_\mu^\perp f_{(1)}) f_{(2)} \ot s_\mu\,,
\end{align}
and fulfils a Leibnitz type formula:
\begin{align}
s_\lambda^\perp(fg) &=\sum_{\mu,\nu} c^\lambda_{\mu,\nu} \,
s_\mu^\perp(f) \, s_\nu^\perp(g)\,,
\end{align}
justifying the name derivative. It is furthermore a rather important fact, that
using the identification $\pi_0=1$, $\pi_n={p_n}\cdot$ and 
$\pi_{-n}=n\partial/\partial_{p_n}$ one easily checks that these operators
generate the {\bf Heisenberg Lie algebra}
\begin{align}\label{Eq-Heisenberg}
  [\pi_n,\pi_m] &= n\delta_{n+m,0}\pi_0 \,,
\end{align} 
closely related to vertex operators and the Witt, and Virasoro algebras used in
string theory. 

The main point we make in this section is to exemplify that in the case of the
character ring Hopf algebras of the classical groups, the notion of 
the \emph{adjoint of multiplication} and that of the \emph{Foulkes
derivative} need no longer be identical; they are logically distinct.
Therefore we need new notation, and we
choose to write $\,^\dagger$ for the adjoint, and keep $\,^\perp$
for the Foulkes derivative.

\subsection{Adjoint of multiplication}
\label{subsec-adj}

\mybenv{Theorem}
The adjoints of multiplication in \CR$\grpGL$, \CR$\grpO$ and \CR$\grpSp$ 
with respect to the general linear, orthogonal and symplectic Schur-Hall scalar 
products are defined to be such that:
\begin{align}
\la \{\nu\} \mid \{\mu\}^\dagger(\{\lambda\}) \ra
&= \la \{\mu\}\cdot \{\nu\}\mid \{\lambda\}\ra \,;\nn
\la [\nu] \mid [\mu]^\dagger([\lambda]) \ra_{2}
&= \la [\mu]\cdot [\nu]\mid [\lambda]\ra_{2} \,;\nn
\la \la\nu\ra \mid \la\mu\ra^\dagger(\la\lambda\ra) \ra_{11}
&= \la \la\mu\ra\cdot \la\nu\ra\mid \la\lambda\ra\ra_{11} \,,
\label{Eq-adj}
\end{align}
respectively, for all partitions $\lambda$, $\mu$ and $\nu$. 
The action of these adjoints then take the explicit forms:
\begin{align}
\{\mu\}^\dagger(\{\lambda\})&=\{\lambda/\mu\}\,; &
[\,\mu\,]^\dagger([\,\lambda\,])&= [\,\mu\,]\cdot[\,\lambda\,]\,; &
\la\mu\ra^\dagger(\la\lambda\ra) &= \la\mu\ra\cdot\la\lambda\ra \,.
\label{Eq-adj-action}
\end{align}
\myeenv

\paragraph{Proof:}
We compute both sides of the requirements (\ref{Eq-adj}) separately
using (\ref{Eq-adj-action}) on the left hand side.
First in the general linear case we have
\begin{align}
 \la \{\nu\} \mid \{\mu\}^\dagger(\{\lambda\})\ra &=\la \{\nu\} \mid \{\lambda/\mu\}\ra
  = \sum_\sigma\ c^\lambda_{\mu,\sigma} \la \{\nu\} \mid \{\sigma\}\ra = c^\lambda_{\mu,\nu}\,,\cr 
  \la \{\mu\}\cdot\{\nu\} \mid \{\lambda\}\ra &= \sum_\rho\ c^\rho_{\mu,\nu} \la \{\rho\}\mid \{\lambda\}\ra = c^\lambda_{\mu,\nu}\,,
\end{align}
so that the two sides are identical as required.
In the orthogonal case we have
\begin{align}
\la [\nu]\mid [\mu]^\dagger([\lambda])\ra_2
&= \sum_{\rho} \la [\nu] \mid [\mu]\cdot[\lambda]\ra_2 
= \sum_{\rho} \la [\nu] \mid [\mu/\rho\cdot \lambda/\rho]\ra_2 \nn 
&= \sum_{\rho,\sigma,\tau,\eta}
   c^\mu_{\rho,\sigma} c^\lambda_{\rho,\tau} c^\eta_{\sigma,\tau}
   \la [\nu]\mid [\eta]\ra_2 
= \sum_{\rho,\sigma,\tau}
   c^\mu_{\rho,\sigma} c^\lambda_{\rho,\tau} c^\nu_{\sigma,\tau}\,,\cr
\la [\mu]\cdot [\nu]\mid [\lambda]\ra_2
&= \sum_{\sigma} \la [\mu/\sigma\cdot \nu/\sigma]\mid [\lambda]\ra_2 \nn
&= \sum_{\sigma,\rho,\tau,\eta}
   c^\mu_{\sigma,\rho} c^\nu_{\sigma,\tau} c^\eta_{\rho,\tau}
   \la [\eta]\mid [\lambda]\ra_2 
= \sum_{\rho,\sigma,\tau}
   c^\mu_{\sigma,\rho} c^\nu_{\sigma,\tau} c^\lambda_{\rho,\tau} \,.
\end{align}
The symmetry $c_{\rho,\sigma}^\mu=c_{\sigma,\rho}^\mu$ then immediately yields 
equality, as required. An entirely analogous proof applies in the symplectic case.
\qed
\paragraph{Remark.}
We are thus left with the fact, that \emph{multiplication is a selfadjoint}
operation in \CR$\grpO$ and in \CR$\grpSp$ with respect to the orthogonal
and symplectic scalar products, respectively. In terms of group representations this
amounts to saying that one can use the second rank tensor $g_{ij} = g_{ji}$
of symmetry  type $\{2\}$ or $f_{ij} = - f_{ji}$ of symmetry type $\{11\}$ to
raise or lower indices. Co- and contra-variant representations of the same index
symmetry type are hence isomorphic.
\qed

\subsection{Foulkes derivative}
\label{Subsec-Foulkes}

To find the correct Foulkes derivative, we exploit both comultiplication and the Schur-Hall 
scalar product in defining any $a^\perp$ as follows
\mybenv{Definition}
The Foulkes derivative is defined in an invariant way as
\begin{align}
a^\perp(b) &= \la a \mid b_{(1)}\ra b_{(2)} \,.
\label{Eq-Foulkes}
\end{align}
\myeenv  
It is easy to check that this definition is equivalent to the skew in the
ordinary $\SY\Lambda$ case. 
\begin{align}
   s_\lambda^\perp(s_\mu) &= \sum_\zeta \la s_\lambda \mid s_\zeta \ra\ s_{\mu/\zeta}
   =\sum_\zeta \delta_{\lambda\zeta}\ s_{\mu/\zeta} = s_{\mu/\lambda}\,.
\end{align}
Furthermore, this definition can be written down in any character Hopf algebra
where we have defined a Schur-Hall scalar product which represents the
orthogonality of irreducible (indecomposable) characters.

\mybenv{Corollary}
The Foulkes derivatives in the case of $\CR\grpGL$, $\CR\grpO$ and $\CR\grpSp$
are given by:
\begin{align}
(s_\lambda)^\perp(s_\mu)
&= \{\lambda\}^\perp(\{\mu\}) =\{\mu/\lambda\}\,;\nn
(o_\lambda)^\perp(o_\mu)
&= \[\lambda\]^\perp(\[\mu\]) = \[\mu/(\lambda\,D)\]\,;\nn
(sp_\lambda)^\perp(sp_\mu) 
&= \la\lambda\ra^\perp(\la\mu\ra) =\la \mu/(\lambda\,B)\ra \,.
\label{Eq-glosp-perp}
\end{align}
\myeenv

\paragraph{Proof:} The Hopf algebra definition for the Foulkes derivative is
basis free, but depends on the scalar product, so that rephrasing (\ref{Eq-Foulkes}) 
in the case of general linear, orthogonal and symplectic characters yields  
\begin{align}
\{\lambda\}^\perp(\{\mu\}) 
&= \la \{\lambda\} \mid \{\mu_{(1)}\} \ra\ \{\mu_{(2)}\}
= \sum_\zeta \la \{\lambda\} \mid \{\zeta\} \ra\ \{\mu/\zeta\}\,,\nn
\[\lambda\]^\perp(\[\mu\]) 
&= \la \[\lambda\] \mid \[\mu_{\[1\]}\]\ra_2 \, \[\mu_{\[2\]}\] 
= \sum_\zeta \la \[\lambda\] \mid \[\zeta\] \ra_2\ \[\mu/(\zeta\,D)\]\,,\nn
\la\lambda\ra^\perp(\la\mu\ra) 
&= \la \la\lambda\ra \mid \la\mu_{\la 1\ra}\ra\ra_{11}\, \la\mu_{\la 2\ra}\ra
= \sum_\zeta \la \la\lambda\ra \mid \la\zeta\ra \ra_{11}\ \la\mu/(\zeta\,B)\ra\,.
\label{Eq-glosp-perp-proof}
\end{align}
where in the case of the orthogonal characters the fourth form of the coproduct 
given in (\ref{Eq-DeltaO-4}) has been used, and its analogue in the case of the 
symplectic characters. 
\qed

It is well known that the above definition (\ref{Eq-Foulkes}) defines a derivation if the
element $a$ is a primitive element in the dual Hopf algebra~\cite{fauser:2002d,fauser:2002c}. Here
applying (\ref{Eq-Foulkes}) in the case of a primitive element ($m\ge 1$) of $\CR\grpO$ we have
\begin{align}
p_n^\perp(p_m) 
&= \la p_n \mid p_m \ra_2 + \la p_n\mid [0]\ra_2 \, p_m 
\end{align}
a consequence of (\ref{Eq-Delta-pn}). However using (\ref{Eq-pn-O}) 
we have 
\begin{align}
\la p_n \mid p_m\ra_2
&= \la \sum_{b=0}^{n-1} (-1)^b [n-b,1^b] + \chi(2\vert n) [0]
   \ \ \ \big\vert\  
           \sum_{d=0}^{m-1} (-1)^d [m-d,1^d] + \chi(2\vert m) [0]  
   \ra_2 \nn
&= n\,\delta_{n,m}\,+\, \chi(2\vert n)\chi(2\vert m)\,,
\end{align}
and
\begin{align}
\la p_n \mid [0]\ra_2
&=
\la \sum_{b=0}^{n-1} (-1)^{b}[n-b,1^b] +\chi(2\vert n)[0] \  \big\vert\
  [0] \ra_2
 = \chi(2\vert n)\,.
\end{align}
Combining these results gives 
\begin{align}
p_n^\perp(p_m) 
&= n\delta_{n,m} + \chi(2\vert m)\chi(2\vert n) + \chi(2\vert n)p_{m}\,.
\label{Eq-pn-perp}
\end{align}

\paragraph{Remark.} This, and an identical result in the symplectic case, 
shows that the power sum basis is not
orthogonal with respect to the orthogonal or symplectic Schur-Hall scalar
products. Furthermore, due to the different Hopf algebra structures of $H$ and
$H^*$, the power sums $p_n$ are {\bf not} the primitive elements of $H^*$. Hence
the identification $p_n^\perp=n \partial/\partial p_n$ of (\ref{Eq-pnperp}) 
that applies in the $\grpGL$ case fails to hold in the $\grpO$ and $\grpSp$ cases. The
correct way to introduce such (formal) derivatives would be to detect the
primitive elements of $H^*$ and to find their dual basis under the relevant
Schur-Hall scalar product. After this identification one could set up 
orthogonal and symplectic Heisenberg Lie algebras quite distinct
from~\eqref{Eq-Heisenberg}. 
This is, however, beyond the scope of the present paper. 
\qed

\subsection{The dual Hopf algebras}
\label{Subsec-dual}

In fact neither \CR$\grpO$ nor \CR$\grpSp$ are 
self-dual Hopf algebras with respect to either the \CR$\grpGL$ 
Schur-Hall scalar product or the \CR$\grpO$, respectively \CR$\grpSp$, 
scalar product. This shows that notwithstanding the Hopf algebra isomorphisms 
between \CR$\grpGL$ and both \CR$\grpO$ and \CR$\grpSp$, these latter
Hopf algebras are not identical to \CR$\grpGL$ since, unlike \CR$\grpGL$ 
they are not self-dual.  Since we will typically consider products such
as $\grpH \ot \grpH^*$ of a Hopf algebra and its dual
(as in the case of the Drinfeld quantum double, or Schur functors with both 
multiplication endomorphisms and Foulkes derivatives, 
or the case of rational characters discussed in section \ref{Sec-rational-universal}),
we note that the branching process does not provide an
isomorphism of this extended structure, and hence the map from one to the other
is a nontrivial transformation.

We now identify convenient bases of the dual Hopf algebras, $\CR\grpO^*$ and 
$\CR\grpSp^*$ of the orthogonal and symplectic character Hopf algebras, $\CR\grpO$ and $\CR\grpSp$, 
respectively, and give explicit formulae for their structure maps. 
Since once more the orthogonal and symplectic cases work out
similarly, we give only the orthogonal versions. Symplectic versions
can be easily obtained by the usual recipe of changing the character brackets
$[~]\rightarrow\la~\ra$ and interchanging series $A\leftrightarrow C$ and
$B\leftrightarrow D$.

\mybenv{Proposition}
Let $\CR\grpO^*$ denote the Hopf algebra dual to $\CR\grpO$. Then a basis
of $\CR\grpO^*$ is provided by the universal characters $[\lambda]^*=\{\lambda\cdot D\}=\{\lambda\}\cdot D$ 
which are such that
\begin{align}
 [\lambda]^*([\mu]) &:= \la [\lambda]^*\mid [\mu] \ra = \delta_{\lambda,\mu}\,.
\end{align}
\myeenv
\paragraph{Proof:}
\begin{align}
   \la [\lambda]^*\mid [\mu] \ra = \la \{\lambda\} D \mid \{\mu/C\} \ra
   = \la \{\lambda\}\cdot DC \mid \{\mu\} \ra = \la \{\lambda\} \mid \{\mu\} \ra = \delta_{\lambda,\mu}\,.
\end{align}
\qed

\mybenv{Proposition}
The dual Hopf algebra \CR$\grpO^*$ is subject to the following structure maps:
\begin{align}
&\textrm{product}
&&&
m([\mu]^*\ot [\nu]^*)&=[\mu]^*\cdot [\nu]^* = [\mu\cdot\nu\cdot D]^* 
\nn\nn[-1ex]{}
&\textrm{unit}
&&&\qquad \eta(1)&= [C]^* \quad\hbox{with}\quad
[C]^*\cdot [\lambda]^* = [\lambda]^* = [\lambda]^*\cdot [C]^* 
 \nn\nn[-1ex]{}
&\textrm{coproduct}
&&&
\delta([\lambda]^*) &= \sum_{\sigma,\zeta} 
  [(\lambda/\sigma)\cdot\zeta]^*  \ot [\sigma\cdot\zeta]^*
\nn\nn[-1ex]{}
&\textrm{counit}
&&&
\epsilon([\lambda]^*) &= \delta_{\lambda,0}
\nn\nn[-1ex]{}
&\textrm{antipode}
&&&
S([\lambda]^*) &= (-1)^{\vert\lambda\vert}[\lambda^\prime BC]^*    
\end{align}
\myeenv
\paragraph{Proof:}
For the {\bf product} we compute
\begin{align}\label{Eq-dual-product}
m([\mu]^*\ot [\nu]^*) = m(\{\mu\cdot D\}\ot\{\nu\cdot D\})=\{\mu\cdot D\cdot\nu\cdot D\}
=\{\mu\cdot\nu\cdot D\}\cdot D= [\mu\cdot\nu\cdot D]^* \,.
\end{align}
For the {\bf unit} we just note that
\begin{align} 
 [C]^*=\{C\cdot D\}=\{0\}\,,
\end{align} 
so that
\begin{align}
[C]^*\cdot[\lambda]^*=\{0\}\cdot\{\lambda\cdot D\}=\{\lambda\cdot D\}=[\lambda]^*
\hbox{~~and~~}
[\lambda]^*\cdot[C]^*=\{\lambda\cdot D\}\cdot\{0\}=\{\lambda\cdot D\}=[\lambda]^*\,.
\end{align}

A little more work is required for the {\bf coproduct}
\begin{align}
\Delta([\lambda]^*)&=\Delta(\{\lambda\cdot D\}) =\Delta(\{\lambda\})\cdot\Delta(D)
=\sum_{\sigma}
  (\{\lambda/\sigma\}\ot\{\sigma\})\cdot(D\ot D)\cdot(\sum_\zeta \{\zeta\}\ot \{\zeta\})\nn
&=\sum_{\zeta,\sigma}
  (\{\lambda/\sigma)\cdot\zeta\}\cdot D)\ot(\{\sigma\cdot\zeta\}\cdot D) 
      =\sum_{\zeta,\sigma}
  [(\lambda/\sigma)\cdot\zeta]^*\ot [\sigma\cdot\zeta]^*\,,
\end{align}
where the coproduct of $D$ has been taken from Proposition~(\ref{Prop-series-coproducts}).

The {\bf counit} maps as follows
\begin{align}
\epsilon([\lambda]^*)=\epsilon(\{\lambda\}\cdot D)=\delta_{\lambda,0}\,.
\end{align}
While the explicit form of the {\bf antipode} action
for the dual Hopf algebra is given by
\begin{align}
S([\lambda]^*)
=S(\{\lambda\cdot D\})
=(-1)^{\vert\lambda\vert} \{\lambda'\cdot D'\}
=(-1)^{\vert\lambda\vert} \{\lambda'\cdot B \cdot CD \}
=(-1)^{\vert\lambda\vert} [\lambda'\cdot BC]^*\,,
\end{align}
where use has been made of the fact that all partitions
in the set $\D=2\P$ are of even weight, and that $D'=B$
\qed
\paragraph{Remark.} A dramatic difference between the character ring Hopf
algebras for orthogonal and symplectic groups and their dual Hopf algebras is
that the product maps of the former are filtered  and hence contain only
finitely many terms, as in (\ref{Eq-NL-O}). 
The dual character ring Hopf algebras, however, have
products based on infinite Schur function series and acquire thereby an
infinite number of terms, as in (\ref{Eq-dual-product}). 
Indeed, the basis elements of these dual Hopf 
algebras, $\[\lambda\]^*=\{\lambda\}\cdot D$ and $\la\lambda\ra^*=\{\lambda\}\cdot B$, 
clearly belong not to the ring $\Lambda$ but to the extension of $\Lambda$ to
include infinite series of Schur functions. In fact these
basis elements are the universal characters
of lowest weight infinite-dimensional holomorphic discrete series
irreducible representations of the $N\rightarrow\infty$ limit of
the non-compact groups $SO^*(N)$ and $Sp(N,\openR)$~\cite{king:wybourne:1985a}.

We add a few (more or less obvious) statements about this structure without
explicit proof.
\mybenv{Corollary}
The dual Hopf algebra \CR$\grpO^*$ is connected, that is we have:
\begin{align}
\Delta(\eta(1))=\Delta([C]^*)=\Delta(\{0\})=\{0\}\ot\{0\}=[C]^*\ot[C]^*
= \eta(1)\ot \eta(1)\,; \nn \nn
\epsilon(m([\lambda]^* \ot [\nu]^*))=\epsilon(\{\lambda\cdot D\cdot \mu\cdot D\})=\delta_{\lambda,0}\,\delta_{\mu,0}
=\epsilon([\lambda]^*)\, \epsilon( [\nu]^*)\,.
\end{align}
\myeenv

However, note that neither the product nor the coproduct is graded. On
the other hand 
the connectedness property allows us to conclude that the antipode still
is an antialgebra homomorphism (though we are bicommutative here), that is
\begin{align}
S(m([\lambda]^*\ot [\mu]^*))&=S(\{\lambda\cdot D\cdot \mu\cdot D\})
=(-1)^{|\lambda|+|\mu|} \{\lambda'\cdot B\cdot \mu'\cdot B\}\nn
&= (-1)^{|\mu|} \{\mu'\cdot BCD\}\ (-1)^{|\lambda|} \{\lambda'\cdot BCD\}\nn
&=(-1)^{|\mu|} [\mu'\cdot BC]^*\ (-1)^{|\lambda|} [\lambda'\cdot BC]^*
=S([\mu]^*)\ S([\lambda]^*)\,,
\end{align} 
where once again use has been made of the fact that all partitions
in the set $\D=2\P$ are of even weight and that $D'=B$.

The fact that the antipode fulfils its defining relation is established
by noting that
\begin{align}
m(1\ot S)\Delta([\lambda]^*)
&=m(1\ot S)\left(\ \sum_{\zeta,\sigma}
   [(\lambda/\sigma)\cdot\zeta]^*\ot [\sigma\cdot\zeta]^*\ \right)\nn
&=m\left(\ \sum_{\zeta,\sigma}
   [(\lambda/\sigma)\cdot\zeta]^*\ot (-1)^{|\sigma|+|\zeta|}[\sigma'\cdot\zeta'\cdot BC]^*\ \right)
 \nn
&=\sum_{\zeta,\sigma} (-1)^{|\sigma|+|\zeta|}
  [(\lambda/\sigma)\cdot\zeta\cdot\sigma'\cdot\zeta'\cdot BC\cdot D]^*
  =\delta_{\lambda,0}\, \sum_\zeta (-1)^{\vert\zeta\vert} [\zeta\cdot\zeta'\cdot B]^* \nn
& =\delta_{\lambda,0}\, [A\,C\,B]^* 
  =\delta_{\lambda,0} [C]^* 
  =\epsilon^*([\lambda]^*)\eta^*(1) \,, 
\end{align}
as required. Use has been made of the antipode 
identity (\ref{Eq-GLai}), $CD=1$, $AB=1$ and the fact that 
$\sum_\zeta (-1)^{\vert\zeta\vert} \{\zeta\cdot\zeta'\}=A\,C$. This
last identity can be established by comparing the dual Cauchy identity~(\ref{Eq-dual-Cauchy})
with the product of the generating functions for the Schur function series $A$ 
and $C$~as given in~(\ref{Eq-series}).

\section{Universal rational characters of the general linear group}
\label{Sec-rational-universal}

There remain further finite-dimensional irreducible representations of these
classical groups. For instance, in the case of $\grpGL(N)$, as well as the 
irreducible covariant tensor representations of highest weight $\lambda$
having character
\begin{align}
   \ch V_{\grpGL(N)}^{\lambda} &=\{\lambda\}(x_1,\ldots,x_N)=s_\mu(x_1,\ldots,x_N)\,,
\label{Eq-charG-covar}
\end{align}
there exist irreducible contravariant tensor 
representations with highest weight $\ov{\mu}=(\ldots,-\mu_2,-\mu_1)$ where
$\mu$ is a partition. These have character
\begin{align}
   \ch V_{\grpGL(N)}^{\ov{\mu}} &=\{\ov{\mu}\}(x_1,\ldots,x_N)=s_\mu(\ov{x}_1,\ldots,\ov{x}_N)\,,
\label{Eq-charG-contra}
\end{align}
with $\ov{x}_i=x_i^{-1}$ for $i=1,2,\ldots,N$. More generally, there exist
irreducible mixed tensor representations of $\grpGL(N)$ with highest weight 
$(\lambda;\ov{\mu})=(\lambda_1,\lambda_2,\ldots,0,\ldots,0,\ldots,-\mu_2,-\mu_1)$ where
$\lambda$ and $\mu$ are both partitions. 
These representations have rational character~\cite{king:1989a,koike:1989a}
\begin{align}
   \ch V_{\grpGL(N)}^{\lambda;\ov{\mu}} &=\{\lambda;\ov{\mu}\}(x_1,\ldots,x_N;x_1,\ldots,x_N) = 
   \sum_{\zeta\in\P} (-1)^{|\zeta|} s_{\lambda/\zeta}(x_1,\ldots,x_N) s_{\mu/\zeta'}(\ov{x}_1,\ldots,\ov{x}_N)\,.
\label{Eq-charG-mixed}
\end{align}

It is straightforward to realise these characters as finite versions of certain universal rational
characters defined in the ring $\Lambda\otimes\ov{\Lambda}$ of symmetric functions with respect to.
Their definition takes the form~\cite{koike:1989a}
\begin{align}
   \{\lambda;\ov{\mu}\} &= (\{\lambda\}\otimes \{\ov{\mu}\})/ \J 
   = \sum_{\zeta\in\F}\ (-1)^{|\zeta|}\ \{\lambda/\zeta\}\otimes \{\ov{\mu/\zeta'}\} 
   \,.
\label{Eq-mixed-universal}
\end{align}
This has as its inverse the identity
\begin{align}
  \{\lambda\}\otimes \{\ov{\mu}\} &= \{\lambda;\ov{\mu}\}/\K = \sum_{\eta\in\F}\ \{\lambda/\eta;\ov{\mu/\eta}\} \,.
\label{Eq-mixed-universal-inverse}
\end{align}

To be more explicit, in terms of two denumerably infinite sequences of indeterminates, 
say $x=(x_1,x_2,\ldots)$ and $y=(y_1,y_2,\ldots)$ we have
\begin{align}
 \{\lambda;\ov{\mu}\}(x;y) &= \sum_{\zeta\in\F}\ (-1)^{|\zeta|}\ s_{\lambda/\zeta}(x)\ s_{\mu/\zeta'}(\ov{y})\,,
\end{align} 
from which we recover our mixed tensor irreducible characters in the form 
\begin{align}
   \ch V_{\grpGL(N)}^{\lambda;\ov{\mu}} &=\{\lambda;\ov{\mu}\}(x_1,\ldots,x_N,0,\ldots,0;x_1,\ldots,x_N,0,\ldots,0) \nn 
 &=\sum_{\zeta\in\P} (-1)^{|\zeta|} s_{\lambda/\zeta}(x_1,\ldots,x_N) s_{\mu/\zeta'}(\ov{x}_1,\ldots,\ov{x}_N)\,,
\label{Eq-charG-mixed-finite}
\end{align}
where we have exploited the usual stability properties of Schur functions with respect
to vanishing indeterminates. 

The notation in (\ref{Eq-mixed-universal}) and (\ref{Eq-mixed-universal-inverse}) is such that in 
$\Lambda\otimes\ov{\Lambda}$ we have:
\begin{align}
  \J=\J_1(x,\ov{y})&=\prod_{i,j}(1-x_i\ov{y_j})=\sum_\zeta (-1)^{|\zeta|}\,s_\zeta(x)\,s_{\zeta'}(\ov{y})
  =\sum_\zeta (-1)^{|\zeta|} \{\zeta\}\ot\{\ov{\zeta'}\}\,;\\ \cr
  \K=\K_1(x,\ov{y})&=\prod_{i,j}(1-x_i\ov{y_j})^{-1}=\sum_\zeta\, s_\zeta(x)\,s_\zeta(\ov{y})
  =\sum_\zeta \{\zeta\}\ot\{\ov{\zeta}\}\,.
\end{align}
where use has been made of the Cauchy identity and its dual. 
Their coproducts take the form
$\Delta(\J)=(\J\otimes\J)\cdot \J'$ and $\Delta(\K)=(\K\otimes\K)\cdot \K'$ with their
cut coproducts given by
\begin{align}
    \J'&=\sum_{\sigma,\tau\in\P} (-1)^{(|\sigma|+|\tau|)}\,  
   (\{\sigma\}\otimes\{\ov{\tau}\})\otimes(\{\tau'\}\otimes\{\ov{\sigma'}\})\,:\\ \cr
    \K'&=\sum_{\sigma,\tau\in\P} \,
   (\{\sigma\}\otimes\{\ov{\tau}\})\otimes(\{\tau\}\otimes\{\ov{\sigma}\})
\end{align}
\label{Eq-Cauchy-coproducts}
This can be seen by taking the product forms of $\J$ and $\K$ and mapping $x$ to $(x,u)$ and $\ov{y}$ to $(\ov{y}\,\ov{v})$. Separating off the products over $x_i\,\ov{y_j}$ and $u_k\,\ov{v_l}$ leaves products
over $x_i\,\ov{v_l}$ and $\ov{y_j}\,u_k$ that can be expanded once again using the Cauchy identity (\ref{Eq-Cauchy}) 
and its dual (\ref{Eq-dual-Cauchy}) to give the required result.

\bigskip

The universal rational characters $\{\lambda;\ov{\mu}\}$ for all partitions $\lambda$ and $\mu$ form a basis
of $\Lambda\otimes\ov{\Lambda}$. Moreover we have  

\mybenv{Theorem}\label{The-GLrat}
The algebra $\CR\grpGLrat$ generated by the universal rational characters $\{\lambda;\ov{\mu}\}$
is a bicommutative Hopf algebra. Its 
structure maps are given by
\begin{align}
&\textrm{product}   & 
m(\{\kappa;\ov{\lambda}\},\{\mu;\ov{\nu}\})&= \{\kappa;\ov{\lambda}\}\cdot\{\mu;\ov{\nu}\} = \sum_{\sigma,\tau\in\P} \{(\kappa/\sigma)\cdot(\mu/\tau);\ov{(\lambda/\tau)\cdot(\nu/\sigma)}\}\nn
&\textrm{unit}      &
\eta(1)&= \{(0);\ov{(0)}\}\nn\nn[-1ex]{}
&\textrm{coproduct} &
\Delta(\{\mu;\ov{\nu}\}) 
&=\sum_{\sigma,\tau,\rho\in\P} \{\mu/\sigma;\ov{\nu/\tau}\}\ot\{\sigma/\rho;\ov{\tau/\rho}\}
\nn
&\textrm{counit}    &
\epsilon(\{\mu;\ov{\nu}\}) &= \delta_{\mu,(0)}\ \delta_{\nu,(0)} \nn\nn[-1ex]{}
&\textrm{antipode}  &
S(\{\mu;\ov{\nu}\}) 
&= (-1)^{\vert\mu\vert+\vert\nu\vert} \sum_\rho \{\mu^\prime/\rho;\ov{\nu^\prime/\rho}\} \nn 
&\textrm{scalar product}  &
\la\cdot\mid\cdot\ra_{1;\ov1}(\{\kappa;\ov{\lambda}\}\ot\{\mu;\ov{\nu}\})
&=\la\{\kappa;\ov{\lambda}\}\,\mid\,\{\mu;\ov{\nu}\}\ra_{1;\ov1} = \delta_{\kappa,\mu}\ \delta_{\lambda,\nu}\,.
\end{align}
\myeenv

The product formula of Theorem~\ref{The-GLrat} was originally given as the rule for decomposing products
of irreducible mixed tensors of $\grpGL$ by Abramsky and King~\cite{abramsky:king:1970a,king:1971a,king:1989a}. 
Its derivation can be accomplished most easily within the
framework of the current paper by noting, precisely as in the derivation 
of the Newell-Littlewood product formula (\ref{Eq-NL-O}),
that $\{\kappa;\ov{\lambda}\}\cdot\{\mu;\ov{\nu}\}
=((\{\kappa\}\otimes\{\ov{\lambda}\})/\J)\cdot((\{\mu\}\otimes\{\ov{\nu}\})/\J)  
=(((\{\kappa\}\otimes\{\ov{\lambda}\})/\J)\cdot((\{\mu\}\otimes\{\ov{\nu}\})/\J)))/\K $ 
where the coefficient of $\{\eta;\ov{\zeta}\}$ in this last expression is given by
\def\nnsp{\nonumber\\[2ex]}
\begin{align}
&\hskip-1cm\la \{\eta\}\otimes\{\ov{\zeta}\} \,\big|\, (((\{\kappa\}\otimes\{\ov{\lambda}\})/\J)\cdot((\{\mu\}\otimes\{\ov{\nu}\})/\J)))/\K  \ra
\nnsp
&= \la (\{\eta\}\otimes\{\ov{\zeta}\}) \cdot K \,|\, ((\{\kappa\}\otimes\{\ov{\lambda}\})/\J)\cdot((\{\mu\}\otimes\{\ov{\nu}\})/\J)  \ra
\nnsp
&= \la \Delta((\{\eta\}\otimes\{\ov{\zeta}\})\cdot K)
\,\big|\, 
((\{\kappa\}\otimes\{\ov{\lambda}\})/\J)\otimes((\{\mu\}\otimes\{\ov{\nu}\})/\J) \ra
\nnsp
&= \la \Delta(\{\eta\}\otimes\{\ov{\zeta}\}) \cdot (K\otimes K)\cdot K' 
\,\big|\, (\{\kappa\}\otimes\{\ov{\lambda}\})/\J)\otimes((\{\mu\}\otimes\{\ov{\nu}\})/\J) \ra
\nnsp
&= \la \Delta(\{\eta\}\otimes\{\ov{\zeta}\})\cdot K'
\,\big|\, ((\{\kappa\}\otimes\{\ov{\lambda}\})/(\J\K))\otimes((\{\mu\}\otimes\{\ov{\nu}\})/(\J\K))  \ra
\nnsp
&= \la \Delta(\{\eta\}\otimes\{\ov{\zeta}\}) \,|\, ( (\{\kappa\}\otimes\{\ov{\lambda}\})\otimes (\{\mu\}\otimes\{\ov{\nu}\}) )/K' \ra
\nnsp
&= \sum_{\sigma,\tau\in\P} \la \Delta(\{\eta\}\otimes\{\ov{\zeta}\}) \,\big|\, 
(\{\kappa/\sigma\}\otimes\{\ov{\lambda/\tau}\})\otimes (\{\mu/\tau\}\otimes\{\ov{\nu/\sigma}\}) \ra
\nn[1ex]
&= \sum_{\sigma,\tau\in\P} \la \{\eta\}\otimes\{\ov{\zeta}\} \,\big|\, 
(\{\kappa/\sigma\}\cdot\{\mu/\tau\})\otimes (\{\ov{\lambda/\tau}\}\cdot\{\ov{\nu/\sigma}\}) \ra
\,,
\end{align}
which gives the product rule of Theorem~\ref{The-GLrat}.

The corresponding result for the evaluation of coproducts coincides 
with the large $M$ and $N$ limit of the
branching rule formula for the restriction from $\grpGL(M+N)$ to $\grpGL(M)\times\grpGL(N)$~\cite{king:1975a}
and may be derived as follows:

\begin{align}
\Delta(\{\mu;\ov{\nu}\} &=  \sum_{\zeta\in\P} (-1)^{|\zeta|} \Delta(\{\mu/\zeta\})\otimes\Delta(\{\ov{\nu/\zeta'}\}) \cr
&= \sum_{\zeta,\sigma,\tau\in\P}(-1)^{|\zeta|} 
(\{\mu/\zeta\sigma\}\otimes\{\sigma\}) \otimes (\{\ov{\nu/\zeta'\tau}\}\otimes\{\ov{\tau}\}) \cr
&= \sum_{\zeta,\sigma,\tau\in\P}(-1)^{|\zeta|} 
(\{\mu/\sigma\zeta\}\otimes\{\ov{\nu/\tau\zeta'}\})) \otimes (\{\sigma\}\otimes\{\ov{\tau}\}) \cr
&= \sum_{\sigma,\tau\in\P} \{\mu/\sigma;\ov{\nu/\tau}\} \otimes \sum_{\rho\in\P} \{\sigma/\rho;\ov{\tau/\rho}\} \,,
\end{align}
which coincides with the required result for the coproduct given in Theorem~\ref{The-GLrat}.

With these two results it is straightforward but rather tedious to verify that the unit, counit and antipode
structure maps displayed in Theorem~\ref{The-GLrat} satisfy all the requirements of a Hopf algebra including the bialgebra and antipode conditions.

\section{Conclusions and discussion}


Our treatment shows that on the Hopf algebraic side the two character
ring Hopf algebras \CR$\grpO$ and \CR$\grpSp$ behave in exactly the same way.
They share the same product structure and differ only in the coproduct where the
series $D$ and $B$ are involved. This stems from the fact that the deformation
of the product actually depends only on the proper cut part $\Delta^\prime$
of the coproduct
\begin{align}
\Delta^\prime(a)
 &= \Delta(a)-1\ot a-a\ot1 \,.
\end{align}
These proper cut parts of $\Delta^\prime(\{2\})$ and
$\Delta^\prime(\{11\})$ are identical (simply the single term $\{1\}\ot\{1\}$), producing the same
deformation. As shown in \cite{fauser:jarvis:king:wybourne:2005a} this is
no longer true for deformations based on tensors of higher degree. For example
$\Delta^\prime(\{3\})$, $\Delta^\prime(\{21\})$ and $\Delta^\prime(\{111\})$ are all different.

The orthogonal and symplectic character of the underlying group finds its
counterpart in the proper definition of the various symmetric function bases.
While the primitives look similar, complete and orthogonal symmetric functions
differ. This is important for applications in physics, since orthogonal,
elementary and power sum symmetric functions can be used to encode partition
functions of physical systems
\cite{schmidt:schnack:2002a,schmidt:schnack:2002b}. Assuming one has a gas of
particles, say atoms or even molecules, having an internal orthogonal or 
symplectic symmetry, one is naturally led to the bases defined in the previous
sections.

We have been able with this approach to evaluate products and coproducts
within the Hopf algebras of the universal character rings of the orthogonal
and symplectic groups, thereby constructing explict formulae for the decomposition 
of products of representations of these groups and of the restriction of these
representations to a variety of subgroups.

To recover the irreducible characters of $\grpO(N)$ and $\grpSp(N)$
in the finite $N$ case one merely limits
the arguments of the universal characters to the eigenvalues of the 
relevant group elements $g$ supplemented by zeros.
Denoting the eigenvalues by $x_k$ and $\xb_k=x_k^{-1}$ for $k=1,2,\ldots,K$,
together with $\pm1$ and $x_{2K+1}$, as appropriate, one obtains:
\begin{align}
\ch V_{\grpO(2K)}^\lambda&=
\[\lambda\](x_1,\ldots,x_K,\xb_1,\ldots,\xb_K,0,\ldots,0)
&&\hbox{for}~~g\in\grpSO(2K)\,;\nn
\ch V_{\grpO(2K)}^\lambda&=
\[\lambda\](x_1,\ldots,x_{K-1},\xb_1,\ldots,\xb_{K-1},1,-1,0,\ldots,0)
&&\hbox{for}~~g\notin\grpSO(2K)\,;\nn
\ch V_{\grpO(2K+1)}^\lambda&=
\[\lambda\](x_1,\ldots,x_K,\xb_1,\ldots,\xb_K,1,0,\ldots,0)
&&\hbox{for}~~g\in\grpSO(2K+1)\,;\nn
\ch V_{\grpO(2K+1)}^\lambda&=
\[\lambda\](x_1,\ldots,x_K,\xb_1,\ldots,\xb_K,-1,0,\ldots,0)
&&\hbox{for}~~g\notin\grpSO(2K+1)\,;\nn
\ch V_{\grpSp(2K)}^\lambda&=
\la\lambda\ra(x_1,\ldots,x_K,\xb_1,\ldots,\xb_K,0,\ldots,0)
&&\hbox{for}~~g\in\grpSp(2K)\,;\nn
\ch V_{\grpSp(2K+1)}^\lambda&=
\la\lambda\ra(x_1,\ldots,x_K,\xb_1,\ldots,\xb_K,x_{2K+1},0,\ldots,0)
&&\hbox{for}~~g\in\grpSp(2K+1)\,,\nn
\label{Eq-chVuniv}
\end{align}
where the final character of $\grpSp(2K+1)$ is indecomposable, rather
than irreducible, with the first $2K$ eigenvalues being those of an
element of $\grpSp(2K)$ and the $x_{2K+1}$ being an element of $\grpGL(1)$.

In order to exploit to the full the results on products and branchings implied by the 
properties of the universal character rings $\CR\grpO$ and $\CR\grpSp$
in the context of the groups $\grpO(N)$ and $\grpSp(N)$ it is necessary to invoke 
certain modification rules~\cite{newell:1951a,king:1971a,black:king:wybourne:1983a,koike:terada:1987a} that 
apply to the above characters whenever the
length $\ell(\lambda)$ of the partitions $\lambda$ exceeds $K$, or in the case
of $SO(2K)$ is equal to $K$.

This is also necessary in dealing with the restriction of the universal
rational characters to the case of $\grpGL(N)$ for finite $N$. As already 
noted, if we denote the eigenvalues of the group element of $\grpGL(N)$ by $x_k$ 
for $k=1,2,\ldots,N$, then the corresponding characters are given by
\begin{align}
   \ch V_{\grpGL(N)}^{\lambda;\ov{\mu}} &=\{\lambda;\ov{\mu}\}(x_1,\ldots,x_N,0,\ldots,0; x_1,\ldots,x_N,0,\ldots,0) \,. 
\label{Eq-mixed-finite}
\end{align}
The corresponding  modification rules have been described elsewhere~\cite{king:1971a,black:king:wybourne:1983a,koike:1989a}.

There remain finite dimensional irreducible spin and indeed spinor
representations of the orthogonal groups that we have not discussed here.
It is possible to make earlier
developments~\cite{king:1975a,black:king:wybourne:1983a} more complete
and rigorous by defining, following Okada~\cite{okada:2010a}, not only
spin characters but also universal spinor characters in the original
ring, $\Lambda$, but now with coefficients in
$\openZ[\varepsilon]/\langle\varepsilon^2-1\rangle$. However, the
double-valued spin characters do not lend themselves in general to a
Hopf algebra analysis because of the complications that arise from the
necessity of distinguishing not only the even and odd $O(2K)$ and
$O(2K+1)$ cases, but also the effect of products, coproducts and 
the antipode map on the spin representations of group elements of
determinant $+1$ and $-1$. In addition, consideration of the characters
of such spin representations of dimension $2^K$ would take us outside
our chosen domain of symmetric functions. 

Applications of symmetric function techniques are widespread. We have
argued in this paper that it is important to pursue the Hopf
algebraic machinery behind the character Hopf algebras, to generalize
these techniques to form a powerful tool which can deal with more
general subsymmetries than orthogonal and symplectic ones. We
restricted our studies here to the classical cases, but even there
novel points arose. 

We believe that the general branching scenario is quite universal,
and have proposed to take it as a blueprint for quantum field
calculations~\cite{fauser:jarvis:2003a,fauser:jarvis:2006a}. The
present work is preparatory to the study of the character ring Hopf
algebras at a (conformal) quantum field level using vertex operator
techniques. It has already allowed the construction of vertex operators for
both the orthogonal and symplectic groups as well as
an extension of such constructions to more general 
subgroups of the general linear group $\grpGL$~\cite{fauser:jarvis:king:2010a}.

Finally we reiterate that the presently developed machinery was
obtained by literally re-doing the quantum field theory calculations
done in~\cite{fauser:2001b,brouder:fauser:frabetti:oeckl:2002a}, in
the context of symmetric functions. We hope to show elsewhere, that
the insights gained here can in turn be profitably applied in quantum
field theory, clarifying algebraic constructions from a group
representation point of view. In particular non-classical subgroup
branchings may lead to new methods allowing the computation of 
nontrivial, that is non-quadratic, invariants, and offering the
possibility of new frameworks for general interacting quantum field
theories.

\section{Acknowledgement}
PDJ and BF acknowledge the Australian Research Council, research
grant DP0208808, for partial support. They also thank the Alexander
von Humboldt  Foundation for a `sur place' travel grant to BF for a
visit to the  University of Tasmania, where this work was done, and
also the School of  Mathematics and Physics for hospitality. RCK is
pleased to acknowledge the award of a Leverhulme Emeritus Fellowship
allowing him to visit the University of Tasmania.
All three authors also gratefully acknowledge a `Research in Pairs' 
grant from the Mathematisches Forschungsinstitut Oberwohlfach, in Spring
2010, which enabled this work to be consolidated, and the warm
hospitality extended to them in Oberwohlfach. 
Finally, this work could not have been completed
without the generous support of the Quantum Computing Group, Department
of Computer Science, University of Oxford, for hosting a research visit.

{\small
%

%
}
\bigskip

\noindent
{\sc
Bertfried Fauser, School of Computer Science, University of Birmingham, \hfil\newline
Edgbaston, Birmingham, B15 2TT, England,
{\small\tt b.fauser@cs.bham.ac.uk}}
\bigskip

\noindent
{\sc
Peter D. Jarvis, University of Tasmania, School of Mathematics and Physics,
Private Bag 37, GPO, Hobart, TAS 7001, Australia, 
{\small\tt Peter.Jarvis@utas.edu.au}}
\bigskip

\noindent
{\sc
Ronald C. King, School of Mathematics, University of Southampton, \hfil\newline
Southampton SO17 1BJ, England, {\small\tt r.c.king@soton.ac.uk}}

\end{document}